\documentclass[a4paper,12pt]{amsart}

\usepackage{amscd, amssymb, amsmath}

\begin{document}

\pagestyle{headings}

\font\fiverm=cmr5 
\input{prepictex}
\input{pictex}
\input{postpictex}

\def\clsp{\overline{\operatorname{span}}}
\def\newspan{{\operatorname{span}}}

\renewcommand{\thefootnote}{\dag}

\newtheorem{thm}{Theorem}[section]
\newtheorem{cor}[thm]{Corollary}
\newtheorem{lem}[thm]{Lemma}
\newtheorem{prop}[thm]{Proposition}
\newtheorem{thm1}{Theorem}

\theoremstyle{definition}
\newtheorem{dfn}[thm]{Definition}
\newtheorem{dfns}[thm]{Definitions}

\theoremstyle{remark}
\newtheorem{rmk}[thm]{Remark}
\newtheorem{rmks}[thm]{Remarks}
\newtheorem{example}[thm]{Example}
\newtheorem{examples}[thm]{Examples}
\newtheorem{note}[thm]{Note}
\newtheorem{notes}[thm]{Notes}

\title{$C^*$-algebras of labelled graphs II - Simplicity results}

\author{Teresa Bates}
\address{Teresa Bates \\
School of Mathematics and Statistics \\
The University of NSW \\
UNSW  Sydney 2052 \\ AUSTRALIA} \email{teresa@unsw.edu.au}

\author{David Pask}
\address{David Pask \\
School of Mathematics and Applied Statistics \\
Austin Keane Building (15) \\
University of Wollongong \\
NSW 2522 \\ AUSTRALIA}
\email{dpask@uow.edu.au}


\thanks{This research was supported by ARC Discovery Project DP0665131, the UNSW Faculty Research Grants Scheme,
and the University of Wollongong}

\begin{abstract}
We prove simplicity and pure infiniteness results for a certain class of labelled graph $C^*$-algebras. We show,
by example, that this class of unital labelled graph $C^*$-algebras is strictly larger than the
class of unital graph $C^*$-algebras.
\end{abstract}

\maketitle

\section{Introduction}


This paper has two main aims. The first is to continue the
development of the $C^*$-algebras of labelled graphs begun in
\cite{bp2} and the second is to provide a tractable example which
illustrates why they are worthy of further study.

A labelled graph is a directed graph $E$ in which the edges have
been labelled by symbols coming from a countable alphabet. By
considering the sequences of labels carried by the bi-infinite paths
in $E$ one obtains a shift space $X$; the labelled graph is then
called a presentation of $X$. A directed graph is a (trivial)
example of a labelled graph, and the shift space it presents is a
shift of finite type (see \cite{lm}). In \cite{bp2} we showed how to
associate a $C^*$-algebra to a labelled space, which consists of a
labelled graph together with a certain collection of subsets of
vertices. By making suitable choices of the labelled spaces it was
shown in \cite[Proposition 5.1, Theorem 6.3]{bp2} that the class of
labelled graph $C^*$-algebras includes graph $C^*$-algebras, the
ultragraph $C^*$-algebras of \cite{t1,t2} and the $C^*$-algebras of
shift spaces in the sense of \cite{m,cm}.

In this paper we shall work almost exclusively with the labelled
spaces which arise in connection with shift spaces. In particular we
shall be interested in identifying key properties of our labelled
spaces which allow us to prove results about the simplicity and pure
infiniteness of the associated $C^*$-algebra (see Theorem
\ref{simple} and Theorem \ref{pureinf}).


Up to now, the examples  of labelled spaces that we have considered
have turned out to have $C^*$-algebras isomorphic to the
$C^*$-algebra of the underlying directed graph (see \cite[Theorem
6.6]{bp2}). In this paper we turn our attention to the question of
whether the class of $C^*$-algebras of labelled spaces that we are
considering is strictly larger than the class of graph
$C^*$-algebras. In section \ref{nogashere} we give presentations of
the Dyck shifts $D_N$ and show that their associated $C^*$-algebras
cannot be unital graph $C^*$-algebras. In section \ref{wannabenogashere} we
present a labelled graph which presents an irreducible non-sofic
shift, whose $C^*$-algebra is simple and purely infinite.


There have now been many papers on the $C^*$-algebras associated to
shift spaces (see \cite{tc,cs,m,m99,ce,cm,bp,bp2} for example). A
drawback to some of the approaches is that the canonical
$C^*$-algebra associated to an irreducible shift space is often not
simple (see \cite[Remark 6.10]{bp2}). We believe that an equally
valid way to study the $C^*$-algebras associated to shift spaces is
to study the $C^*$-algebras of the labelled graphs which present
them. This belief is founded on the observation that the labelled
graph $( E_1 , \mathcal{L}_1 )$ of Examples \ref{agreex} (i) is a
presentation of an irreducible sofic shift (called the even shift)
whose $C^*$-algebra is simple (see \cite[Remark 6.10]{bp2}) whereas
the $C^*$-algebra associated to the even shift in \cite{cm} is not
simple.


The work of Matsumoto on symbolic matrix systems and their
associated $\lambda$-graph systems gives us an important method for
studying shift spaces using labelled graphs (see \cite{ma2, ma3,
m99, m100} amongst others). However, we feel that there is an extra
facility afforded by our approach. Whilst $\lambda$-graph systems
are indeed labelled graphs, they are quite complicated. This makes
them difficult to visualise; for instance the labelled graphs in
Examples \ref{agreex} (i) give rise to the same $C^*$-algebras as
the ones for the symbolic matrix systems described on
\cite[p.297]{m99}. Furthermore we believe that our presentations of
the Dyck shifts in section \ref{nogashere} give us a more tractable
way of studying them. Of equal importance is the fact that our
labelled spaces are ideally suited to handle shift spaces over
countably infinite alphabets.


The paper begins with a long section in which we describe many of
the important concepts associated to labelled graphs and labelled
spaces. The two main results of this section are Proposition
\ref{wasitfinitenoughforyou} and Proposition \ref{essmatsys}. In
Proposition \ref{wasitfinitenoughforyou} we give an important
embellishment to the treatment of labelled spaces in \cite{bp2} by
identifying the basic objects in a labelled space, which we call the
generalised vertices. In Proposition \ref{essmatsys} we establish
concrete connections between our work and that of Matsumoto by
showing how to associate a symbolic matrix system to a labelled
graph.

In section \ref{CoLS} we recall the definition of the $C^*$-algebra
of a labelled space from \cite{bp2}. In Proposition \ref{newspandef}
we give a new description of the canonical spanning set for a
labelled graph $C^*$-algebra in terms of generalised vertices. Then
in Proposition \ref{wegotmat} we use this new description to show
the relationship between the $C^*$-algebra of a labelled graph and
the $\lambda$-graph $C^*$-algebra of the associated symbolic matrix
system.

In section \ref{core} we give a description of the AF core of a
labelled graph $C^*$-algebra before moving on to prove the
Cuntz-Krieger uniqueness Theorem (Theorem \ref{ckunique}) in section
\ref{ckut}. The central hypothesis to the Cuntz-Krieger uniqueness
Theorem for labelled graphs is the notion of disagreeability, which
replaces the aperiodicity hypothesis in the corresponding theorem
for directed graphs (see, for example \cite[Theorem 3.1]{bprsz}).

In section \ref{biggies} we give the simplicity and pure
infiniteness results for labelled graph $C^*$-algebras. To prove the
simplicity result (Theorem \ref{simple}) we need a notion of
cofinality appropriate for labelled graphs. The notion of cofinality
for labelled graphs is much more subtle than that for directed
graphs as many different infinite paths in the underlying directed
graph can have the same labels. To prove the pure infiniteness
result (Theorem \ref{pureinf}) we need to examine how periodic paths
arise labelled graphs. The situtation is much more complicated than
for directed graphs since periodic points in the shift space
associated to a labelled graph need not arise from a loop in the
underlying directed graph.

Finally in section \ref{newex} we provide two new examples of
labelled graphs to which our main results apply. In section
\ref{nogashere} we provide a labelled graph presentation of the Dyck
shifts $D_N$. In Proposition \ref{dyckprop} show that these
presentations give rise to simple purely infinite labelled graph
$C^*$-algebras. In Remark \ref{matsumotoagain} we give a formula for
the K-theory of our labelled graph $C^*$-algebras which demonstrates
that the $C^*$-algebras we associate to Dyck shifts cannot be
isomorphic to graph $C^*$-algebras. In section
\ref{wannabenogashere} we provide a presentation of an interesting
new irreducible non-sofic shift whose labelled graph $C^*$-algebra
is simple and purely infinite.

\subsection*{Acknowledgements} The authors would like to acknowledge the
hospitality given to us by the University of Victoria in Canada, the CRM
in Barcelona  and the Fields Institute in Toronto during the preparation of
this paper.

\section{Collected definitions and notation}

\subsection*{Directed graphs} A directed graph $E$ consists of a quadruple $( E^0 , E^1
, r , s )$ where $E^0$ and $E^1$ are (not necessarily countable) sets of vertices and
edges respectively and $r, s : E^1 \to E^0$ are maps giving the
direction of each edge. A path $\lambda = e_1 \ldots e_n$ is a
sequence of edges $e_i \in E^1$ such that $r ( e_i ) = s ( e_{i+1}
)$ for $i=1 , \ldots , n-1$, we define $s ( \lambda ) = s ( e_1 )$ and
$r ( \lambda ) = r ( e_n )$. The collection of paths of length $n$
in $E$ is denoted $E^n$ and the collection of all finite paths in
$E$ by $E^*$, so that $E^* = \bigcup_{n \ge 0} E^n$.

A {\em loop} in $E$ is a path which begins and ends at the same vertex, that is
$\lambda \in E^*$ with $s ( \lambda ) =  r ( \lambda )$. We say that $E$ is {\em row-finite}
if every vertex emits finitely many edges. The graph $E$ is called {\em transitive} if there
given any pair of vertices $u,v \in E^0$ there is a path $\lambda \in E^*$ with
$s(\lambda)=u$ and $r ( \lambda ) = v$. We denote the
collection of all infinite paths in $E$ by $E^\infty$.

\subsection*{Standing assumption 1} We will assume
that our directed graphs $E$ are {\it essential}: all vertices
emit and receive edges (i.e.\ $E$ has
no sinks or sources).

\subsection*{Labelled graphs} A {\it labelled graph} $( E , \mathcal{L}  )$ over a countable alphabet
${\mathcal A}$ consists of a directed graph $E$ together with a
labelling map $\mathcal{L}  : E^1 \to \mathcal{A}$. Without loss of generality we may assume that the map
$\mathcal{L} $ is onto.

Let $\mathcal{A}^*$ be the collection of all {\em words} in the
symbols of $\mathcal{A}$. The map $\mathcal{L} $ extends naturally
to a map $\mathcal{L}  : E^n \to \mathcal{A}^*$, where $n \ge 1$:
for $\lambda = e_1 \ldots e_n \in E^n$ put $\mathcal{L}  ( \lambda )
= \mathcal{L}  ( e_1 ) \ldots \mathcal{L} ( e_n )$; in this case the
path $\lambda \in E^n$ is said to be a {\em representative} of the
{\em labelled path} $\mathcal{L}  ( e_1 ) \ldots \mathcal{L}  ( e_n
)$. Let $\mathcal{L} ( E^n )$ denote the collection of all labelled
paths in $(E,\mathcal{L} )$ of length $n$ where we write $|\alpha|=
n$ if $\alpha \in \mathcal{L}( E^n )$. The set $\mathcal{L}^* (E ) =
\bigcup_{n \ge 1} \mathcal{L} ( E^n )$ is the collection of all
labelled paths in the labelled graph $(E, \mathcal{L}  )$.  We may
similarly extend $\mathcal{L}$ to $E^\infty$.

The labelled graph $( E , \mathcal{L}  )$ is {\em left-resolving} if
for all $v \in E^0$ the map $\mathcal{L}  : r^{-1} (v) \to
\mathcal{A}$ is injective. The left-resolving condition ensures that for all $v \in
E^0$ the labels $\{ \mathcal{L}  (e) : r(e) = v \}$ of all incoming
edges to $v$ are all different. For $\alpha$ in $\mathcal{L}^* ( E )$ we put
\[
s_\mathcal{L}  ( \alpha ) = \{ s ( \lambda ) \in E^0 : \mathcal{L} (
\lambda ) = \alpha \} \text{ and } r_\mathcal{L}  ( \alpha ) = \{ r
( \lambda ) \in E^0 : \mathcal{L}  ( \lambda ) = \alpha \} ,
\]

\noindent so that $r_\mathcal{L}  , s_\mathcal{L}  : \mathcal{L}^* (
E ) \to 2^{E^0}$. We shall drop the subscript on $r_\mathcal{L} $
and $s_\mathcal{L} $ if the context in which it is being used is
clear.

Let $( E, \mathcal{L} )$ be a labelled graph. For $A \subseteq E^0$
and $\alpha \in \mathcal{L}^* (E)$ the {\em relative range of
$\alpha$ with respect to $A$} is defined to be
\[
r_{\mathcal L}(A,\alpha) = \{ r ( \lambda ) : \lambda \in E^* , \mathcal{L}  (
\lambda ) = \alpha , s ( \lambda ) \in A \} .
\]

\noindent A collection $\mathcal{B}  \subseteq 2^{E^0}$ of subsets
of $E^0$ is said to be {\em closed under relative ranges for
$(E , \mathcal{L} )$} if for all $A \in \mathcal{B} $ and $\alpha \in
\mathcal{L}^* (E )$ we have $r (A,\alpha) \in \mathcal{B} $. If
$\mathcal{B} $ is closed under relative ranges for $(E , \mathcal{L}
)$, contains $r ( \alpha )$ for all $\alpha \in \mathcal{L}^* (E  )$
and is also closed under finite intersections and unions, then we
say that $\mathcal{B} $ is {\em accommodating} for $(E, \mathcal{L}
)$.

Let $\mathcal{E}^{0,-}$ denote the smallest subset of $2^{E^0}$
which is accommodating for $( E,\mathcal{L} )$. Since $\mathcal{E}^{0,-}$ is generated
by a countable family of subsets of $E^0$, under countable operations, it follows that $\mathcal{E}^{0,-}$
is countable, even though $E^0$ itself may be uncountable. Of course, $2^{E^0}$ is the largest
accommodating collection of subsets  for $(E,{\mathcal L})$.

\subsection*{Labelled spaces} A {\em labelled space} consists of a triple $(E , \mathcal{L}  ,
\mathcal{B}  )$, where $(E , \mathcal{L}  )$ is a labelled graph and $\mathcal{B} $
is accommodating for $(E , \mathcal{L}  )$ .

A labelled space $(E , \mathcal{L}  , \mathcal{B}  )$ is {\em weakly left-resolving}
if for every $A , B \in \mathcal{B} $ and every $\alpha \in
\mathcal{L}^* ( E )$ we have $r ( A , \alpha ) \cap r ( B , \alpha )
= r ( A \cap B , \alpha)$.

\begin{rmks} \label{reducereuserecycle}
\begin{itemize}
\item[(i)] If  $(E , \mathcal{L} , {\mathcal E}^{0,-})$ is weakly left-resolving then $\mathcal{E}^{0,-}$
is the closure of $\{r(\alpha)\;:\; \alpha \in {\mathcal L}^*(E)\}$ under finite unions and
intersections (cf.\ \cite[Remark 3.9]{bp2}).  Moreover, if $\alpha \in {\mathcal L}^*(E)$ and $A =
\bigcup_{k=1}^m \bigcap_{i=1}^n r ( \beta_{i,k} ) \in
\mathcal{E}^{0,-}$ where $\beta_{i,k} \in \mathcal{L}^* (E)$ for
$i=1 , \ldots , n$ and $k=1 , \ldots , m$ then we have
\[
r ( A , \alpha ) = \bigcup_{k=1}^m \bigcap_{i=1}^n r ( \beta_{i,k}
\alpha ) .
\]

\item[(ii)] If $( E , \mathcal{L}  )$ is left-resolving then the labelled space $( E
, \mathcal{L}  , \mathcal{B}  )$ is weakly left-resolving for any
$\mathcal{B} $.

\end{itemize}
\end{rmks}

For $\ell \ge 1$ and $A \subseteq E^0$ let $E^\ell A = \{ \lambda
\in E^\ell : r ( \lambda ) \in A \}$. The labelled space $( E ,
\mathcal{L} , \mathcal{B} )$ is {\em receiver set-finite} if for all
$A \in \mathcal{B} $ and all $\ell \ge 1$ the set $\mathcal{L} (
E^\ell A ) := \{{\mathcal L}(\lambda) \;:\; \lambda \in E^\ell A \}$
is finite. In particular, the labelled space $( E , \mathcal{L} ,
\mathcal{B} )$ is receiver set-finite if each $A \in {\mathcal B}$
receives only finitely labelled paths of length $\ell$ (even though
$A$ may receive infinitely many paths of each length $\ell$). More
generally, for $\ell \ge 1$ and $A \subseteq E^0$ let
\[
\mathcal{L} ( E^{\le \ell} ) = \bigcup_{j=1}^\ell \mathcal{L} ( E^j
A ) = \bigcup_{j=1}^\ell \mathcal{L} ( E^j  ) \text{ and }
\mathcal{L} ( E^{\le \ell} A ) .
\]

\noindent For $A \subseteq E^0$ and $n \ge 1$ we define $L^n_A = \{
{\mathcal L}(\lambda) \;:\; \lambda \in E^n, s(\lambda) \in A \}$.
If $L^1_A$ is finite for all $A \in {\mathcal B}$  we say that $( E
, \mathcal{L}  , \mathcal{B}  )$ is {\em set-finite}.

\subsection*{Standing assumption 2} We will assume that $(E,{\mathcal
L},{\mathcal E}^{0,-})$ is receiver set-finite, set-finite and
weakly left-resolving.

\begin{rmk}
The conditions of set-finiteness and receiver set-finiteness are
trivially satisfied by labelled spaces over finite alphabets.  The
condition of set-finiteness for labelled spaces is the analogue of
row-finiteness for directed graphs.  Taken together the conditions
of set-finiteness and receiver set-finiteness give us the analogue
of local finiteness for directed graphs.
\end{rmk}
\vspace{3mm}
\noindent

\subsection*{Generalised vertices} For $v \in E^0$ and $\ell \ge 1$ let
\[
\Lambda_\ell (v) = \{ \lambda \in \mathcal{L} ( E^{\le \ell} ) : v
\in r ( \lambda ) \} = \mathcal{L} ( E^{\le \ell} v ) .
\]
The relation $\thicksim_\ell$ on $E^0$ is defined by $v
\thicksim_\ell w$ if and only if $\Lambda_\ell (v) = \Lambda_\ell
(w)$; hence $v \thicksim_\ell w$ if $v$ and $w$ receive exactly the
same labelled paths of length at most $\ell$. Evidently
$\thicksim_\ell$ is an equivalence relation and we use $[v]_\ell$ to
denote the equivalence class of $v \in E^0$. We call the $[v]_\ell$
{\em generalised vertices} as they play the same role in labelled
spaces as vertices in a directed graph.

Set $\Omega_\ell = E^0 / \thicksim_\ell$ and $\Omega := \bigcup_{\ell \ge 1} \Omega_\ell$.
If the alphabet $\mathcal{A}$ is finite, then $\Omega_\ell$ is finite.
If there is $L \ge 1$ such that $\Omega_\ell = \Omega_L$ for all $\ell \ge L$,
then the underlying shift $\textsf{X}_{E,\mathcal{L}}$ is a sofic
shift (see \cite{cm},\cite{lm}). Conversely, if $\textsf{X}$ is a sofic shift then
every presentation $( E , \mathcal{L} )$ of the shift $\textsf{X}$
has this property (see \cite[Exercise 3.2.6]{lm}).

For $\ell \ge 1$ let ${\mathcal E}^{0,-}_\ell \subseteq {\mathcal
E}^{0,-}$ be the smallest  subset of $2^{E^0}$ which
contains $r(\lambda)$ for all $\lambda \in \mathcal{L}(E^{\le
\ell})$ and is closed under finite intersections and unions.
Evidently $\mathcal{E}^{0,-}_\ell \subseteq
\mathcal{E}^{0,-}_{\ell + 1}$.  Following Remark
\ref{reducereuserecycle}(i) we have $\mathcal{E}^{0,-} = \bigcup_{\ell
= 1}^\infty \mathcal{E}^{0,-}_\ell$.

For $v \in E^0$ and $\ell \ge 1$, the equivalence class $[ v ]_\ell$
does not necessarily belong to $\mathcal{E}^{0,-}_\ell$; however, as
we shall see in Proposition \ref{wasitfinitenoughforyou} (i) $[ v
]_\ell$ may be expressed as a difference of elements of
$\mathcal{E}^{0,-}_\ell$. First we need the following technical
lemma.

\begin{lem} \label{isitfinitenoughforyou}
Let $(E,{\mathcal L},{\mathcal E}^{0,-})$ be a labelled space, $v
\in E^0$ and $\ell \ge 1$.
\begin{itemize}
\item[(i)] The set $\Lambda_\ell(v)$ is finite and  $X_\ell(v) :=
\bigcap_{\lambda \in \Lambda_\ell (v)} r ( \lambda )
\in {\mathcal E}^{0,-}_\ell$.  Moreover $[ v ]_\ell \subseteq X_\ell
(v)$.

\item[(ii)]The set of labels
$Y_\ell (v) := \bigcup_{w \in X_\ell (v)} \Lambda_\ell (w) \backslash
\Lambda_\ell (v)$ is finite, and $r ( Y_\ell (v ) ) \in \mathcal{E}^{0,-}_\ell$.
\end{itemize}
\end{lem}

\begin{proof}
For the first statement let $A \in \mathcal{E}^{0,-}$ be such that
$v \in A$. Since $( E , \mathcal{L} , \mathcal{E}^{0,-} )$ is
receiver set-finite $\mathcal{L} ( E^j v ) \subseteq \mathcal{L} ( E
^j A )$ is finite for all $j \ge 1$ and hence $\Lambda_\ell (v)=
\bigcup_{j=1}^\ell {\mathcal L}(E^jv)$ is finite for all $\ell \ge
1$. It now follows that $X_\ell (v)$ is a finite intersection of
elements of $\mathcal{E}^{0,-}_\ell$ and hence $X_\ell (v) \in
\mathcal{E}^{0,-}_\ell$.  Since $X_\ell (v)$ is the set of vertices
which receive at least the same labelled paths as $v$ up to length
$\ell$ we certainly have $[v]_\ell \subseteq X_\ell (v)$.

For the second statement observe that $Y_\ell (v) = \mathcal{L}(
E^{\le \ell} X_\ell (v) ) \backslash \Lambda_\ell (v)$.  Since
$(E,{\mathcal L},{\mathcal E}^{0,-})$ is receiver set-finite and
$X_\ell(v) \in {\mathcal E}^{0,-}$ the sets $\mathcal{L}(E^{\le
\ell}X_\ell(v))$ and $Y_\ell (v)$ must be finite.
Note that $ r ( Y_\ell (v) )  = \bigcup_{\mu \in Y_\ell (v)} r ( \mu )$ belongs
to $\mathcal{E}^{0,-}_\ell$ as it is a finite union of
elements of $\mathcal{E}^{0,-}_\ell$.
\end{proof}

\noindent  The set $Y_\ell (v)$ denotes the additional labelled
paths of length at most $\ell$ received by those vertices which
receive at least the same labelled paths as $v$ up to length $\ell$.

\begin{prop} \label{wasitfinitenoughforyou}
Let $(E,{\mathcal L},{\mathcal E}^{0,-})$ be a  labelled space, $v
\in E^0$ and $\ell \ge 1$.
\begin{itemize}
\item[(i)] We have
$[v]_\ell =  X_\ell (v) \backslash r ( Y_\ell (v) )$.
\item[(ii)]  For every set $A \in {\mathcal E}^{0,-}_\ell$ we can
find vertices $v_1,\dots,v_m \in A$ such that $A = \bigcup_{i=1}^m
[v_i]_\ell$.
\item[(iii)] There are $w_1 , \ldots , w_m \in [v]_\ell$ such that
$[v]_\ell = \bigcup_{i=1}^m [ w_i ]_{\ell + 1}$.
\end{itemize}
\end{prop}

\begin{proof}
For the first statement  observe that $[ v ]_\ell$ consists of those
vertices which  receive exactly the labelled paths from $\Lambda_\ell (v)$
whereas other vertices in $X_\ell (v)$ may receive more labelled
paths. Hence, to form $[v]_\ell$ we  remove those vertices from
$X_\ell( v)$ which receive different labelled paths of length $\ell$ from $v$ --
these are precisely the vertices in $r ( Y_\ell (v) )$.

For the second statement note that by Remark
\ref{reducereuserecycle}(i) any $A \in \mathcal{E}^{0,-}_\ell$ can
be written as a finite union of elements of the form $B_k =
\bigcap_{i=1}^n r ( \beta_i )$ where $\beta_i \in \mathcal{L} (
E^{\le \ell} )$. If $v_1 \in B_k$ then $[ v_1 ]_\ell \subseteq B_k$
as $v_1$, and hence every vertex in $[ v_1 ]_\ell$, must receive
$\beta_1 , \ldots , \beta_n$ and so lie in $B_k$. If $B_k \neq [ v_1
]_\ell$, there is $v_2 \in B_k$ with $\Lambda_\ell (v_1 ) \neq
\Lambda_\ell ( v_2 )$. Again we have $[v_2]_\ell \subseteq B_k$.
Since $(E,{\mathcal L},{\mathcal E}^{0,-})$ is receiver set-finite
$B_k \in \mathcal{E}^{0,-}_\ell$ receives only finitely many
different labelled paths of length at most $\ell$. Hence there are
vertices $\{ v_i : 1 \le i \le m \}$ in $B_k$ such that $B_k =
\bigcup_{i=1}^m [ v_i ]_\ell$ and our result is established.

For the final statement we observe that since
$\mathcal{E}^{0,-}_\ell \subseteq \mathcal{E}^{0,-}_{\ell+1}$ the
first statement  shows that $[v]_\ell$ may be written as a
difference $A \backslash B$ of elements of
$\mathcal{E}^{0,-}_{\ell+1}$. The result then follows by applying
the second statement to $A,B \in \mathcal{E}^{0,-}_{\ell+1}$ and
noting that the $[w_i]_{\ell +1}$'s are disjoint.
\end{proof}

\subsection*{Shift spaces} Let $(E, \mathcal{L} )$  be a labelled graph.
The subshift $\textsf{X}_E$ is defined by $\textsf{X}_E = \{ x \in
(E^1)^{\bf Z} \;:\; s(x_{i+1}) = r(x_i) \mbox{ for all } i \in {\bf
Z } \}$. The subshift $( \textsf{X}_{E,\mathcal{L} } , \sigma )$ is
defined by
\[
\textsf{X}_{E,\mathcal{L}} = \{ y \in \mathcal{A}^{\bf Z} :
\text{ there exists
 } x \in \textsf{X}_E \mbox{ such that } y_i = \mathcal{L}  ( x_i ) \text{ for all } i \in {\bf Z}
\} ,
\]

\noindent where $\sigma$ is the shift map $\sigma(y)_i = y_{i+1}$
for $i \in {\bf Z}$. The labelled graph $(E ,\mathcal{L}  )$ is said
to be a {\em presentation} of the shift space $\textsf{X}_{ E ,
\mathcal{L}  }$ with language ${\mathcal L}^*(E)$.

We are primarily interested in one-sided shift spaces, namely
\[
\textsf{X}_{E,\mathcal{L}}^+= \{ y \in \mathcal{A}^{\bf N} :
\text{ there exists
 } x \in E^\infty \mbox{ such that } y_i = \mathcal{L}  ( x_i ) \text{ for all } i \in {\bf N}
\}
\]

\noindent and we restrict the shift map to $\textsf{X}_{
E,\mathcal{L} }^+$.  For an infinite labelled path $x \in \textsf{X}_{E,\mathcal{L}}^+$
we define $s_{\mathcal{L}}(x)$ to be the set of all $v \in E^0$ for
which there is an infinite path $\widehat{x} \in E^\infty$ with
$s(\widehat{x}) = v$ and ${\mathcal L}(\widehat{x}) = x$. The infinite
path $\hat{x}$ is said to be a {\em representative} of $x$.

An infinite labelled path $x \in \textsf{X}^+_{E,\mathcal{L}}$ is {\em
periodic} if $\sigma^n x = x$ for some $n \ge 1$. A path which is
not periodic is called {\em aperiodic}.

\begin{example} \label{trivialex} If $E$ is a directed graph then we may consider it
as a labelled graph when endowed with the trivial labelling
$\mathcal{L}_t$. In this case $\mathcal{E}^{0,-}$ consists of all
finite subsets of $E^0$ (see \cite[Examples 4.3(i)]{bp2}) and
$[v]_\ell = \{ v \}$ for all $\ell \ge 1$. We shall identify
$\mathcal{L}_t^* ( E )$ with $E^*$ and
$\textsf{X}^+_{E,\mathcal{L}_t}$ with $E^\infty$.
\end{example}

\subsection*{Symbolic Matrix Systems} Essential
symbolic matrix systems are defined in \cite[\S 2]{m99}. To a
left-resolving labelled graph $( E , \mathcal{L} )$ over a finite
alphabet we associate matrices $( M(E)_{\ell,\ell+1} ,
I(E)_{\ell,\ell+1} )_{\ell \ge 1}$  as follows: For $\ell \ge 1$,
write $\Omega_\ell = \{ [v_i]_\ell : i = 1 , \ldots , m ( \ell )
\}$, then $I(E)_{\ell,\ell+1}$ is a $m (\ell ) \times m ( \ell + 1)$
matrix with entries $0,1$ determined by
\begin{equation} \label{iedef}
I(E)_{\ell,\ell+1} ( [v_i]_\ell , [w_j]_{\ell+1} ) =
\begin{cases} 1 & \text{ if } [w_j]_{\ell+1} \subset [v_i]_\ell \\
0 & \text{ otherwise. }
\end{cases}
\end{equation}

\noindent The symbolic matrix $M(E)_{\ell,\ell+1}$ is the same size
as $I(E)_{\ell,\ell+1}$ with entries determined as follows: For $v
\in E^0$ let $\langle v \rangle_\ell$ denote the collection of
labelled paths of length exactly $\ell$ which arrive at $v$. Since
$(E,\mathcal{L})$ is left-resolving we may partition the set of
labelled paths of length $\ell+1$ arriving at $w$ to write $\langle
w \rangle_{\ell+1}$ as the disjoint union
\[
\langle w \rangle_{\ell+1} =  \bigcup_{e \in r^{-1} ( w ) } \langle s(e)
\rangle_\ell \mathcal{L} (e) ,
\]

\noindent where $\langle s(e) \rangle_\ell \mathcal{L} (e)$ denotes
the set of labelled paths of length $\ell+1$ formed by the
juxtaposition of the symbol $\mathcal{L} (e)$ at the end of each
labelled path in $\langle s(e) \rangle_\ell$. Since all vertices in
$[v_i]_\ell$ and $[w_j]_{\ell+1}$ receive the same labelled paths of
length $\ell$ and $\ell+1$ respectively we may unambiguously define
\begin{equation} \label{medef}
M(E)_{\ell,\ell+1} ( [v_i]_\ell , [w_j]_{\ell+1} ) = \sum_{e \in s^{-1} (v_i)\cap r^{-1} (w_j)} \mathcal{L} (e)
\end{equation}

\noindent where the right hand-side is treated as a formal sum.

\begin{prop} \label{essmatsys}
Let $( E , \mathcal{L} )$ be a left-resolving labelled graph over a finite alphabet. Then the matrices
$( M(E)_{\ell,\ell+1} , I(E)_{\ell,\ell+1} )_{\ell \ge 1}$  defined above form an essential symbolic matrix
system.
\end{prop}

\begin{proof}
If suffices to check that the matrices $( M(E)_{\ell,\ell+1} ,
I(E)_{\ell,\ell+1} )_{\ell \ge 1}$ satisfy the conditions on
\cite[p.290]{m99}: Since $E$ is essential it is straightforward to
check from the definition of $I(E)_{\ell,\ell + 1}$ and
$M(E)_{\ell,\ell + 1}$ that conditions (1), (2), (2-a), (2-b), (3),
(5-i) and (5-ii) are satisfied. It remains to check  that for $\ell
\ge 1$ we have $M(E)_{\ell,\ell+1} I(E)_{\ell+1,\ell+2} =
I(E)_{\ell,\ell+1} M(E)_{\ell+1,\ell+2}$.

For $\ell \ge 1$ we form the entry $M(E)_{\ell,\ell+1}
I(E)_{\ell+1,\ell+2} ( [u_i]_\ell , [w_k]_{\ell+1} )$ as follows:
For each $[v_j]_{\ell+1}$ which receives an edge from $[u_i]_\ell$,
the entry is the formal sum of the labels received by the unique
$[v_j]_{\ell+1}$ of which $[w_k]_{\ell+2}$ is a subset. In which
case
\[
M(E)_{\ell,\ell+1} I(E)_{\ell+1,\ell+2} ( [u_i]_\ell , [w_k]_{\ell+2} ) = \sum_{e \in s^{-1} (u_i) \cap r^{-1} (w_k)} \mathcal{L} (e) .
\]

\noindent On the other hand, to form the entry $I(E)_{\ell,\ell+1} M(E)_{\ell+1,\ell+2} ( [u_i]_\ell , [w_k]_{\ell+2} )$ we take
each $[v_j]_{\ell+1}$ which is a subset of $[u_i]_\ell$ and then formally sum the labels of the edges to $[w_k]_{\ell+1}$. In which case
\begin{eqnarray*}
I(E)_{\ell,\ell+1} M(E)_{\ell+1,\ell+2} ( [u_i]_\ell , [w_k]_{\ell+2} )
&= \sum_{[v_j]_{\ell+1} \subseteq [u_i]_{\ell}} \sum_{e \in s^{-1} (v_j) \cap r^{-1} (w_k)} \mathcal{L} (e) \\
&= \sum_{e \in s^{-1} (u_i) \cap r^{-1} (w_k)} \mathcal{L} (e) .
\end{eqnarray*}

\noindent Hence for $\ell \ge 1$ we have $M(E)_{\ell,\ell+1} I(E)_{\ell+1,\ell+2} = I(E)_{\ell,\ell+1} M(E)_{\ell+1,\ell+2}$
as required.
\end{proof}

\section{$C^*$-algebras of labelled spaces} \label{CoLS}

\noindent We recall from \cite{bp2} the definition of the definition
of the $C^*$-algebra associated to the labelled space
$(E,\mathcal{L},\mathcal{E}^{0,-})$.

\begin{dfn} \label{lgdef}
Let $( E , \mathcal{L}  , \mathcal{E}^{0,-}  )$ be a labelled space.
A representation of $( E , \mathcal{L} , \mathcal{E}^{0,-} )$
consists of projections $\{ p_A : A \in \mathcal{E}^{0,-} \}$ and
partial isometries $\{ s_a : a \in \mathcal{A} \}$ with the
properties that

\begin{itemize}

\item[(i)] If $A, B \in \mathcal{E}^{0,-} $ then $p_A p_B = p_{A \cap
B}$ and $p_{A \cup B} = p_A + p_B - p_{A \cap B}$, where
$p_\emptyset = 0$.

\item[(ii)] If $a \in \mathcal{A}$ and $A \in
\mathcal{E}^{0,-} $ then $p_A s_a = s_a p_{r ( A, a )}$.

\item[(iii)] If $a , b \in \mathcal{A}$ then $s_a^*
s_a = p_{r( a )}$ and $s_{a}^*s_{b} = 0$ unless $a = b$.

\item[(iv)] For $A \in \mathcal{E}^{0,-} $ we have
\begin{equation} \label{sumcond}
p_A = \sum_{a \in L^1_A} s_{a} p_{r( A , a )} s_{a}^* .
\end{equation}
\end{itemize}
\end{dfn}

\begin{rmk}
If the directed graph $E$ contains sinks then
we need to modify condition (iv) above (note that ${\mathcal E}^{0,-}$
is different in this case (see \cite[Definition 3.8]{bp2})).   The original definition
\cite[Definition 4.1]{bp2} was in error since it would lead to degeneracy of the
vertex projections for sinks.  We thank Toke Carlsen for pointing this out
to us.  If $A$ contains a finite number of sinks and ${\mathcal B}= {\mathcal E}^{0,-}$
or ${\mathcal E}^0$, then we obtain
the relation
\[
p_A = \sum_{a \in L^1_A} s_{a} p_{r( A , a )} s_{a}^* + \sum_{v \in A\; :\; v\mbox{ is a sink }} p_v.
\]
\end{rmk}

\begin{dfn}
Let $( E , \mathcal{L}  , \mathcal{E}^{0,-}  )$ be a  labelled
space, then $C^* ( E , \mathcal{L}  , \mathcal{E}^{0,-}  )$ is the
universal $C^*$-algebra generated by a representation of $( E ,
\mathcal{L} , \mathcal{E}^{0,-}  )$.
\end{dfn}

\noindent The universal property of $C^*(E, \mathcal{L},
\mathcal{E}^{0,-})$ allows us to define a strongly continuous action
$\gamma$ of ${\bf T}$ on $C^*(E,\mathcal{L},\mathcal{E}^{0,-})$
called the {\em gauge action} (see \cite[Section 5]{bp2}).  As in
\cite[Proposition 3.2]{cbms} we denote by $\Phi$ the conditional
expectation of $C^*( E,\mathcal{L},\mathcal{E}^{0,-})$ onto the
fixed point algebra $C^*(E,\mathcal{L},\mathcal{E}^{0,-})^\gamma$.
 If $(E,\mathcal{L},\mathcal{E}^{0,-})$ is a labelled space then by
\cite[Lemma 4.4]{bp2} we have
\[
C^*(E,{\mathcal L},{\mathcal E}^{0,-}) =
\overline{\text{span}}\{s_\alpha p_A s_\beta^* \;:\; \alpha,\beta
\in {\mathcal L}^*(E), A \in {\mathcal E}^{0,-}\}.
\]

\noindent Indeed, we can write down a more informative spanning set for
$C^*(E,{\mathcal L},{\mathcal E}^{0,-} )$.

\begin{prop} \label{newspandef}
Let $( E , \mathcal{L} , \mathcal{E}^{0,-})$ be a labelled space.
Then
\[
C^* ( E , \mathcal{L} , \mathcal{E}^{0,-} ) = \overline{\text{span}}\{s_\alpha p_{[v]_\ell} s_\beta^* \;:\;
\alpha, \beta \in {\mathcal L}^*(E), [v]_\ell \in \Omega_\ell \}
\]

\noindent where
\begin{equation} \label{relationforvell}
p_{[v]_\ell} := p_{X_\ell (v)} - p_{r ( Y_\ell (v) )} p_{X_\ell (v)}
= \sum_{a \in L^1_{[v]_\ell}} s_a p_{r([v]_\ell,a)}s_a^* .
\end{equation}
\end{prop}

\begin{proof}
The first assertion holds from repeated applications of Proposition
\ref{wasitfinitenoughforyou}. Applying \eqref{sumcond} of Definition
\ref{lgdef} we have
\[
p_{[v]_\ell} = p_{X_\ell(v)} - p_{X_\ell(v) \cap r(Y_\ell(v))}
             = \sum_{a \in L^1_{X_\ell(v)}} s_a
             p_{r(X_\ell(v),a)}s_a^* - \sum_{b \in L^1_{X_\ell(v)
             \cap r(Y_\ell(v))}} s_b p_{r(X_\ell(v)
             \cap r(Y_\ell(v)),b)} s_b^*.
\]
In order to eliminate double counting of labels that are emitted by
both $X_\ell(v)$ and $r(Y_\ell(v))$ we need to split
$L^1_{X_\ell(v)}$ into two disjoint parts (the labels that come only
out of $X_\ell(v)$ and those that come out of both $X_\ell(v)$ and
$r(Y_\ell(v))$) to obtain
\[
p_{[v]_\ell}= \sum_{a \in L^1_{X_\ell(v)} \backslash L^1_{X_\ell(v)
\backslash r(Y_\ell(v))}}
             \!\!\! s_a p_{r(X_\ell(v),a)}s_a^* + \sum_{b \in L^1_{X_\ell(v)
             \cap r(Y_\ell(v))}} s_b \left(p_{r(X_\ell(v),b)} - p_{r(X_\ell(v)
             \cap r(Y_\ell(v)),b)}\right) s_b^*.
\]
We may replace $X_\ell( v)$ in the first sum by $[v]_\ell$ as the
labels $a$ are emitted only by the vertices in $[v]_\ell$ and not by
the vertices in $X_\ell(v) \cap r(Y_\ell(v))$.  In the second sum
the labels $b$ are emitted by both $[v]_\ell$ and $X_\ell(v) \cap
r(Y_\ell(v))$, but we subtract the projections corresponding to the
copies emitted by $X_\ell(v) \cap r(Y_\ell(v))$ and so we have
equation \eqref{relationforvell} as required.
\end{proof}

\begin{rmk}
Note that while proving Proposition \ref{newspandef} we have shown
that for $[v]_\ell \in \Omega_\ell$ and $a \in {\mathcal A}$
\[
r([v]_\ell,a) = r(X_\ell(v),a) \backslash r(r(Y_\ell(v)),a)
\]
which is a difference of two elements of ${\mathcal E}^{0,-}_{\ell +
1}$.
\end{rmk}
\noindent Recall from Proposition \ref{essmatsys} that to a labelled
graph $( E , \mathcal{L} )$ over a finite alphabet $\mathcal{A}$ we
may associate an essential symbolic matrix system $(
M(E)_{\ell,\ell+1} , I(E)_{\ell,\ell+1} )_{\ell \ge 1}$. By
\cite[Proposition 2.1]{m99} there is a unique (up to isomorphism)
$\lambda$-graph system $\mathfrak{L}_{E,\mathcal{L}}$ associated to
$( M(E)_{\ell,\ell+1} , I(E)_{\ell,\ell+1} )_{\ell \ge 1}$. By
\cite[Theorem 3.6]{m100} one may associate a $C^*$-algebra
$\mathcal{O}_{\mathfrak{L}_{E,\mathcal{L}}}$ to the $\lambda$-graph
system $\mathfrak{L}_{E,\mathcal{L}}$ which is the universal
$C^*$-algebra generated by partial isometries $\{ t_a : a \in
\mathcal{A} \}$ and projections $\{ E_i^\ell : i = 1 , \ldots , m (
\ell ) \}$ satisfying relations
\begin{align}
\sum_{a \in \mathcal{A}} t_a t_a^*  = 1 & \hphantom{abcdefghijklmnop} \label{idrel} \\
\sum_{i=1}^{m(\ell)} E_i^\ell = 1 & \hphantom{abcdef}
E_i^\ell = \sum_{j=1}^{m(\ell+1)} I (E)_{\ell,\ell+1} (i,j) E_j^{\ell+1} \text{ for } i = 1 , \ldots , m ( \ell ) \label{vellrel} \\
t_a t_a^* E_i^\ell &= E_i^\ell t_a t_a^* \text{ for } a \in {\mathcal A} \text{ and } i = 1 , \ldots , m ( \ell ) \label{commute} \\
t_a^* E_i^\ell t_a &= \sum_{j=1}^{m(\ell+1)} A_{\ell , \ell+1} ( i ,
a , j ) E_j^{\ell+1}  \text{ for } a \in {\mathcal A} \text{ and } i
= 1 , \ldots , m ( \ell ) \label{ckrel}
\end{align}

\noindent where $A_{\ell,\ell+1} (i,a,j)=1$ if $a$ occurs in the formal sum $M(E)_{\ell,\ell+1} ( [v_i ]_\ell , [v_j]_{\ell+1} )$
and is $0$ otherwise.

\begin{prop} \label{wegotmat}
Let $(E , \mathcal{L} )$ be a left-resolving labelled graph over a
finite alphabet. Then we have $C^* ( E , \mathcal{L} ,
\mathcal{E}^{0,-} ) \cong
\mathcal{O}_{\mathfrak{L}_{E,\mathcal{L}}}$ where
$\mathfrak{L}_{E,\mathcal{L}}$ is the $\lambda$-graph system
associated to the symbolic matrix system $( M(E)_{\ell,\ell+1} ,
I(E)_{\ell,\ell+1} )_{\ell \ge 1}$.
\end{prop}

\begin{proof}
By Proposition \ref{newspandef} the elements $\{s_a \;:\; a \in
{\mathcal A}\}$ and $\{p_{[v_i]_\ell} \;:\; i = 1, \dots, m(\ell)\}$
form a generating set for $C^*(E,{\mathcal L},{\mathcal E}^{0,-})$.
Let $T_a = s_a$ and $F^\ell_i = p_{[v_i]_\ell}$ then $\{T_a,
F^\ell_i\}$ satisfy relations \eqref{idrel} -- \eqref{ckrel} above.
Hence by the universal property of
$\mathcal{O}_{\mathfrak{L}_{E,\mathcal{L}}}$   there is a map
$\pi_{T,F} : \mathcal{O}_{\mathfrak{L}_{E,\mathcal{L}}} \to
C^*(E,{\mathcal L},{\mathcal E}^{0,-})$ characterised by $\pi_{T,F}
(t_a) = T_a$ and $\pi_{T,F}(E^\ell_i) = F^\ell_i$.

Let $\{t_a : a \in {\mathcal A}\}$ and $\{E^\ell_i \;:\; i =
1,\dots,m(\ell)\}$ be generators for
$\mathcal{O}_{\mathfrak{L}_{E,\mathcal{L}}}$.   For $A \in {\mathcal
E}^{0,-}_\ell$ and $a \in {\mathcal A}$ let $P_A =
\sum_{i\;:\;[v_i]_\ell \subseteq A} E_i^\ell$ and $S_a = t_a$.  One
checks that $\{S_a,P_A\}$ is a representation of the labelled space
$(E,\mathcal{L},{\mathcal E}^{0,-})$. By universality of
$C^*(E,\mathcal{L},\mathcal{E}^{0,-})$ there is a map $\pi_{S,P} :
C^* ( E , \mathcal{L} , \mathcal{E}^{0,-} ) \to
\mathcal{O}_{\mathfrak{L}_{E,\mathcal{L}}}$ characterised by
$\pi_{S,P}(s_a) = S_a$ and $\pi_{S,P}(p_A) = P_A$. In particular, we
have $\pi_{S,P}(p_{[v_i]_\ell}) = P_{[v_i]_\ell}$ for all $i \in
1,\dots,m(\ell)$.  Our result follows since $\pi_{T,F}$ and
$\pi_{S,P}$ are inverses of one another.

\end{proof}

\section{AF core} \label{core}

\noindent
In this section we perform a detailed analysis of the AF core of $C^* ( E , \mathcal{L} , \mathcal{E}^{0,-} )$
which we will need to prove the main result of the following section.

\begin{dfn}
\noindent For $1 \le k \le \ell$ let
\[
\mathcal{F}^k ( \ell ) = \overline{\text{span}} \{ s_\alpha p_A
s_\beta^* : \alpha , \beta \in \mathcal{L} ( E^k ) , A \in
\mathcal{E}^{0,-}_\ell \}.
\]
\end{dfn}

\noindent For $\ell \ge 1$ and $[v]_\ell \in \Omega_\ell $ we have
$p_{[v]_\ell} \in \mathcal{F}^k ( \ell )$ as $X_\ell (v) , r (
Y_\ell (v) ) \in \mathcal{E}^{0,-}_\ell$ by Lemma
\ref{isitfinitenoughforyou} (ii).

\begin{dfn}
For $1 \le k \le \ell$ and $[v]_\ell \in \Omega_\ell$ let
\[
\mathcal{F}^k ( [v]_\ell ) = \overline{\text{span}} \{ s_\alpha
p_{[v]_\ell} s_\beta^* : \alpha , \beta \in \mathcal{L} ( E^k ) \}.
\]
\end{dfn}

\begin{prop} \label{include}
For $1 \le k \le \ell$ we have

\begin{itemize}
\item[(i)] $\mathcal{F}^k ( \ell ) \cong \oplus_{[v]_\ell} \mathcal{F}^k ( [v]_\ell )$,
where each $\mathcal{F}^k ( [v]_\ell )$ is a finite-dimensional
matrix algebra.
\item[(ii)] For each  $v \in E^0$ there  are $w_1,\dots,w_n \in
[v]_\ell$ such that ${\mathcal F}^k([v]_\ell) = \oplus_{i=1}^n
{\mathcal F}^k([w_i]_{\ell + 1})$.  Hence $\mathcal{F}^k ( \ell )
\subseteq \mathcal{F}^k ( \ell+1 )$.
\item[(iii)] There
is an embedding of $\mathcal{F}^k ( \ell )$ into $\mathcal{F}^{k+1}
( \ell+1 )$.
\end{itemize}
\end{prop}

\begin{proof}
For the first statement of (i), applying Proposition
\ref{wasitfinitenoughforyou} (ii) shows that every element $s_\alpha
p_A s_\beta^* \in \mathcal{F}^k ( \ell )$ can be written as a finite
sum of elements of the form $s_\alpha p_{[v]_\ell} s_\beta^* \in
\mathcal{F}^k ( [v]_\ell )$. The result follows as the summands in
the decomposition are mutually orthogonal since $|\alpha| = |\beta|
= k$ and the equivalence classes $[v]_\ell$ are disjoint. For the
second statement of (i) note that since $[v]_\ell$ can be written as
the difference of two elements of ${\mathcal E}^{0,-}$ it receives
only finitely many different labelled paths of length $k$ and hence
the set $\{s_\alpha p_{[v]_\ell} s_\beta^*\;:\; |\alpha| = |\beta| =
k\}$ is finite. It is straightforward to show that the elements
$s_\alpha p_{[v]_\ell} s_\beta^*$ form a system of matrix units in
$\mathcal{F}^k ( [v]_\ell )$ and the result follows.

Part (ii) follows by Proposition \ref{wasitfinitenoughforyou} (iii).
Part (iii) follows from Definition \ref{lgdef} (iv).
\end{proof}

\begin{thm}
Let $(E , \mathcal{L} , \mathcal{E}^{0,-} )$ be a labelled space,
then $\mathcal{F} = \overline{ \bigcup_{k,\ell}
\mathcal{F}^k(\ell)}$ is an AF algebra with $\mathcal{F} \cong C^*
(E , \mathcal{L} , \mathcal{E}^{0,-} )^\gamma$.
\end{thm}

\begin{proof}
The first statement follows from Proposition \ref{include}.  The
second statement follows by an argument similar to that of
\cite[Lemma 2.2]{bprsz}.
\end{proof}

\section{Cuntz-Krieger Uniqueness Theorem} \label{ckut}

Recall from \cite[\S 3]{kpr} that the directed graph $E$ satisfies
condition (L) if every loop has an exit; that is if $\lambda \in
E^n$ is a loop, then there is some $1 \le i \le n$ such that the
vertex $r ( \lambda_i )$ emits more than one edge. Condition (L) is
the key hypothesis for the Cuntz-Krieger uniqueness theorem for
directed graphs (see \cite[Theorem 3.7]{kpr}, \cite[Theorem
3.1]{bprsz}). Since periodic paths in $E^\infty$ arise from loops in
$E$, condition (L) guarantees that there are lots of paths in
$E^\infty$ which aperiodic.

In this section we seek an analogue for condition (L) in the context
of labelled graphs which will allow us to prove a Cuntz-Krieger
uniqueness theorem for labelled graph $C^*$-algebras.  The correct
analogue for condition (L) must ensure the existence of aperiodic
paths in $\textsf{X}^+_{E,{\mathcal L}}$.   The two key difficulties
to overcome in the context of labelled graphs are that we must
accommodate the generalised vertices $[v]_\ell$ in a labelled graph
and deal with the fact that a periodic path $x \in
\textsf{X}^+_{E,\mathcal{L}}$ need not arise from a loop in $E$.

The following definition is inspired by \cite[Lemma 3.7]{cbms}.

\begin{dfns}
Let $( E , \mathcal{L} , \mathcal{E}^{0,-})$ be a labelled space,
 $[v]_\ell \in \Omega_\ell$ and  $\alpha \in {\mathcal
L}^*(E)$ be such that $|\alpha|>1$ and $s(\alpha) \cap [v]_\ell \ne
\emptyset$.  We say that $\alpha$ is {\em agreeable} for $[v]_\ell$
if there are $\alpha' , \beta, \gamma \in {\mathcal L}^*(E)$ with $|
\beta | = | \gamma | \le \ell$ and $\alpha = \beta \alpha' = \alpha'
\gamma$. Otherwise we say that $\alpha$ is {\em disagreeable} for
$[v]_\ell$.

We say that $[v]_\ell$ is {\em disagreeable} if there is an $N
>0$ such that for all $n > N$ there is an $\alpha \in {\mathcal
L}^*(E)$ with $\vert \alpha \vert \ge n$ that is disagreeable for
$[v]_\ell$.

The labelled space $(E,\mathcal{L},\mathcal{E}^{0,-})$ is {\em
disagreeable} if for every $v \in E^0$ there is an $L_v > 0$ such
that $[v]_\ell$ is disagreeable  for all $\ell > L_v$.
\end{dfns}

\begin{rmk} \label{wlog}
Suppose that $[v]_p = \bigcup_{i=1}^m [w_i]_q$, where $q > p$.  Then
$L_v \ge L_{w_i}$ for all $i \in \{1,\dots,m\}$ since each $w_i \in
[v]_p$.
\end{rmk}

\noindent The following Lemma shows that the notion of disagreeability reduces to
condition (L) for directed graphs and so is the appropriate
condition for us to use in our Cuntz-Krieger uniqueness theorem and
simplicity results.

\begin{lem}
The directed graph $E$ satisfies condition (L) if and only if the
labelled space $(E, \mathcal{L}_t , \mathcal{E}^{0,-})$ is
disagreeable.
\end{lem}

\begin{proof}
Suppose that $E$ satisfies condition (L). Observe that for all $\ell
\ge 1$ and all $v \in E^0$, $[v]_{\ell} = \{v\}$.  We show that
every $v \in E^0$ is disagreeable.  Let $L_v = 1$, $N=1$, fix $n >
N$ and $\ell > L_v$.  If $v$ does not lie on a  loop, then any path
$\alpha$ with $|\alpha| \ge n$ is disagreeable for $[v]_\ell =
\{v\}$.  If $v$ does lie on a loop $\alpha = \alpha_1 \dots
\alpha_m$, without loss of generality we may assume that $s(\alpha)
= v$.  Since $E$ satisfies condition $(L)$ there is a path $\beta$
with $s(\beta) = v$ and $\beta_{|\beta|} \not \in
\{\alpha_1,\dots,\alpha_m\}$.  The path $\alpha^n \beta$ has length
$\ge n$ and is disagreeable for $[v]_\ell$.

Suppose $E$ does not satisfy condition (L).  Then there is a $v \in
E^0$ and a simple loop $\alpha$ with $s(\alpha) = v$ that has no
exit.  Let $N > 0$.  Then there is an $n$ such that $|\alpha^n|
> N$. Suppose $n \ge 2$.  We claim that $\lambda = \alpha^n$ is
agreeable for every $\ell > |\alpha|$.  Set $\beta = \gamma =
\alpha$ and $\lambda' = \alpha^{n-1}$.  Since  $\lambda = \beta \lambda'
= \lambda' \gamma$ where $|\beta| = |\gamma| \le \ell$  it follows that
$[v]_\ell = \{ v \}$ is agreeable
for $\ell$.  Since $\alpha^n$ is the only path of length $n|\alpha|$
emitted by $v$, it follows that $v$ is not disagreeable.  Thus the
labelled space $(E, \mathcal{L}, \mathcal{E}^{0,-})$ is not
disagreeable.
\end{proof}

\begin{examples} \label{agreex}
\begin{itemize}
\item[(i)] Recall from \cite[Examples 3.3 (iii)]{bp2} the labelled graphs
\begin{equation*}
\beginpicture
\setcoordinatesystem units <1cm,0.75cm>

\setplotarea x from -1 to 12, y from -0.75 to 1.1

\put{$(E_1 , \mathcal{L} _1 ) :=$}[l] at -0.75 0

\put{$\bullet$} at 3 0

\put{$\bullet$} at 5 0

\put{$1$}[l] at 1.6 0

\put{$0$}[b] at 4 1.1

\put{$0$}[t] at 4 -1.1

\put{$u$}[l] at 3.15 0

\put{$v$}[l] at 5.15 0

\setquadratic

\plot 3.1 0.1  4 1 4.9 0.1 /

\plot 3.1 -0.1 4 -1  4.9 -0.1 /

\circulararc 360 degrees from 3 0 center at 2.5 0

\arrow <0.25cm> [0.1,0.3] from 2.015 0.1 to 2 -0.1 

\arrow <0.25cm> [0.1,0.3] from 4.1 -0.985 to 3.9 -1

\arrow <0.25cm> [0.1,0.3] from 3.9 0.985 to 4.1 1

\put{$(E_2 , \mathcal{L} _2 ) :=$} at 7 0

\put{$\bullet$} at 10 0

\put{$\bullet$} at 12 0

\put{$\bullet$} at 10 -2

\put{$1$}[l] at 8.6 0

\put{$0$}[b] at 11 1.1

\put{$0$}[t] at 11 -1.1

\put{$u$}[l] at 10.15 0

\put{$v$}[l] at 12.15 0

\put{$w$}[bl] at 10.15 -1.95

\put{$1$}[r] at 9.85 -1

\put{$0$}[r] at 9.47 -2.5

\setquadratic

\plot 10.1 0.1  11 1 11.9 0.1 /

\plot 10.1 -0.1 11 -1  11.9 -0.1 /

\circulararc 360 degrees from 10 0 center at 9.5 0

\circulararc 360 degrees from 10 -2 center at 10 -2.5

\arrow <0.25cm> [0.1,0.3] from 9.015 0.1 to 9 -0.1 

\arrow <0.25cm> [0.1,0.3] from 11.1 -0.985 to 10.9 -1

\arrow <0.25cm> [0.1,0.3] from 10.9 0.985 to 11.1 1

\arrow <0.25cm> [0.1,0.3] from 10 -0.15 to 10 -1.85

\arrow <0.25cm> [0.1,0.3] from 9.95 -2.98 to 10.04 -3.01
\endpicture
\end{equation*}

\noindent
are set-finite, receiver set-finite, left-resolving presentations of
the even shift.

Consider $(E_1,\mathcal{L}_1)$.  We claim that
$(E_1,\mathcal{L}_1,\mathcal{E}^{0,-}_1)$ is disagreeable.  Now for
all $\ell \ge 1$ we have $[u]_\ell = \{u\}$. Let $L_u = 1$ and $N =3
$.  Then for $n > N$ the labelled path $\alpha_n = 1 1^n 0$
satisfies $|\alpha_n| = n + 2 \ge N$ and $\alpha_n$ is disagreeable
for $[u]_\ell$ as its first and last symbols disagree. Also for all
$\ell \ge 1$ we have $[v]_\ell = \{v\}$. If we let $N = 4$ and $L_v
= 1$, then for each $n > N$ the path $\alpha_n = 0^{2n+1}1$
satisfies  $|\alpha_n| = 2n + 2 \ge n$ and $\alpha_n$ is
disagreeable for $[v]_\ell$ as its first and last symbols disagree.
Thus the labelled space
 $(E_1,\mathcal{L}_1,\mathcal{E}^{0,-}_1)$ is disagreeable and our claim is
established.

Consider $(E_2,\mathcal{L}_2)$.  We claim that $[w]_\ell$ is
agreeable for all $\ell \ge 2$. Now  $[w]_\ell = \{w\}$ for all
$\ell \ge 2$, and any labelled path $\alpha$ satisfying $s(\alpha)
\cap [w]_\ell \ne \emptyset$ must have the form $\alpha = 0^n$ for
some $n$.  But $\alpha = 0^n$ is agreeable for $[w]_\ell$ for all
$\ell \ge 2$ whenever $n \ge \ell + 1$: set $\alpha' = 0^{n-\ell}$,
$\beta = \gamma = 0$. Thus $(E_2,\mathcal{L}_2,\mathcal{E}^{0,-}_2)$
is not disagreeable.

\item[(ii)] Let $G$ be a group with a finite set of generators
$S = \{ g_1 , \ldots , g_m \}$, such that $g_i \neq g_j$ for $i \neq
j$. The (right) {\it Cayley graph of $G$ with respect to $S$} is the
essential row-finite directed graph $E_{G,S}$ where $E_{G,S}^0 = G$,
$E_{G,S}^1 = G \times  S$ with range and source maps given by $r (h,
g_i) = hg_i$ and $s (h, g_i) = h$ for $i=1, \ldots , m$. The map
$\mathcal{L}_{G,S} (h,g_i)=g_i$ gives us a set-finite, receiver
set-finite, labelled graph $( E_{G,S} , \mathcal{L}_{G,S} )$. Since
$G$ is cancelative it follows that $( E_{G,S} , \mathcal{L}_{G,S} )$
is left resolving. As each vertex in $E_{G,S}$ receives the same
labelled paths it follows that $[g]_\ell = G$ for all $g \in G$ and
$\ell \ge 1$ and so ${\mathcal E}_{G,S}^{0,-} = \{ \emptyset, G \}$.
Each $g \in G$ emits the same $m^\ell$ labelled paths of length
$\ell$. So if $m = \vert S \vert > 1$, it follows that for all
$[g]_\ell = G$ there is a disagreeable labelled path of length $n>1$
beginning at $[g]_\ell =G$. Hence $( E_{G,S}, \mathcal{L}_{G,S},
\mathcal{E}_{G,S}^{0,-} )$ is disagreeable.
\end{itemize}
\end{examples}

\begin{thm}\label{ckunique}
Let $(E , \mathcal{L} , \mathcal{E}^{0,-} )$ be a labelled space. If
$\{ T_\alpha , Q_A \}$ and $\{ S_\alpha , P_A \}$ are two
representations of $(E , \mathcal{L} , \mathcal{E}^{0,-} )$ in which
all the projections $p_A , P_A$ are nonzero, then there is an
isomorphism $\phi$ of $C^* ( T_\alpha , Q_A )$ onto $C^* ( S_\alpha
, P_A )$ such that $\phi ( T_\alpha ) = S_\alpha$ and $\phi ( Q_A )
= P_A$.
\end{thm}

\noindent To prove this theorem we show that the representations $\pi_{T,Q}$
and $\pi_{S,P}$ of $C^*(E,{\mathcal L},{\mathcal E}^{0,-})$ are
faithful. The required isomorphism will then be $\phi = \pi_{S,P}
\circ \pi_{ T, Q}^{-1}$.  The usual approach is to invoke symmetry
and prove that
\begin{itemize}
\item[(a)] $\pi_{S,P}$ is faithful on $C^*(E,{\mathcal
L},{\mathcal E}^{0,-})^\gamma$ and
\item[(b)] $\|\pi_{S,P} \left( \Phi (a) \right)\|
\le \|\pi_{S,P}(a)\|$ for all $a \in C^*(E,{\mathcal L},{\mathcal
E}^{0,-})$.
\end{itemize}

\noindent Part (a) is proved in \cite[Theorem 5.3]{bp2}.  To prove
(b) we must do a little more work than is needed for graph
$C^*$-algebras because of the more complicated structure of
$C^*(E,{\mathcal L},{\mathcal E}^{0,-})^\gamma$ as is discussed in
section \ref{core}.

\begin{proof}
By Proposition \ref{newspandef} every element of $C^*(E,{\mathcal
L},{\mathcal E}^{0,-})$ may be approximated by elements of the form
\[
a= \sum_{(\alpha , [w]_\ell , \beta ) \in F} c_{\alpha , [w]_\ell ,
\beta} s_\alpha p_{[w]_\ell} s_\beta^*
\]

\noindent where $F$ is finite, and so it is enough to prove (b) for
such elements $a$.

Let $k = \max \{ \vert \alpha \vert , \vert \beta \vert : ( \alpha ,
[w]_\ell , \beta ) \in F\}$ .  By Proposition \ref{newspandef} we
may suppose (changing $F$ if necessary), that every $(\alpha,
[w]_\ell , \beta) \in F$ is such that $\min\{\vert \alpha \vert,
\vert \beta \vert : ( \alpha , [w]_\ell , \beta ) \in F \} = k$. Let
$M = \max\{\vert \alpha \vert, \vert \beta \vert: (\alpha , [w]_\ell
, \beta) \in F\}$
and $L = \max \{ L_w : (\alpha,[w]_\ell, \beta) \in F \}$.  By
Remark \ref{wlog} and Proposition \ref{wasitfinitenoughforyou}(iii)
we may suppose (again changing $F$ if necessary, but not $M$ or $k$)
that $\ell \ge \max\{L,M-k\}$.

Since $\vert \alpha \vert = \vert \beta \vert$ implies that $\vert
\alpha \vert = k$ we have
\[
\Phi (a) = \sum_{(\alpha , [w]_\ell , \beta ) \in F , \vert \alpha
\vert = \vert \beta \vert} c_{\alpha , [w]_\ell , \beta} s_\alpha
p_{[w]_\ell} s_\beta^* \in \mathcal{F}^k ( \ell )
\]
\noindent where $\Phi$ is the conditional expectation of
$C^*(E,\mathcal{L},\mathcal{E}^{0,-})$ onto
$C^*(E,\mathcal{L},\mathcal{E}^{0,-})^\gamma$. By Proposition
\ref{include} (i) $\mathcal{F}^k ( \ell )$ decomposes as the
$C^*$-algebraic direct sum $\oplus_{[w]_\ell} \mathcal{F}^k (
[w]_\ell )$, so does its image under $\pi_{S,P}$, and there is a
$[v]_\ell \in \Omega_\ell$ such that $\|\pi_{S,P}(\Phi(a))\|$ is
attained on ${\mathcal F}^k([v]_\ell)$.  Let $F_{[v]_\ell}$ denote
the elements of $F$ of the form  $(\alpha,[v]_\ell,\beta)$, then we
have
\[
\| \pi_{S,P} ( \Phi (a) ) \| = \left\| \sum_{(\alpha , [v]_\ell ,
\beta ) \in F_{[v]_\ell} , \vert \alpha \vert = \vert \beta \vert}
\!\!\!\! c_{\alpha,[v]_\ell,\beta} S_\alpha P_{[v]_\ell} S_\beta^*
\right\|.
\]

\noindent We write
\[
b_v = \sum_{(\alpha , [v]_\ell , \beta ) \in F_{[v]_\ell} , \vert
\alpha \vert = \vert \beta \vert} c_{\alpha , [v]_\ell , \beta}
S_\alpha P_{[v]_\ell} S_\beta^*
\]

\noindent
and let $G = \{ \alpha : \text{ either } (\alpha,[v]_\ell,\beta) \in
F_{[v]_\ell} \text{ or } (\beta,[v]_\ell,\alpha) \in F_{[v]_\ell}
\text{ with } |\alpha| = |\beta| \}$. Then $\operatorname{span} \{
S_\alpha P_{[v]_\ell} S_\beta^* : \alpha , \beta  \in G \}$ is a
finite dimensional matrix algebra containing $b_v$.

Since $\ell > L$,  $[ v ]_\ell$  is disagreeable. Hence there is an
$n > M$ and a $\lambda$ with $\vert \lambda \vert \ge n$  such that
$[v]_\ell \cap s ( \lambda ) \ne \emptyset$ which has no
factorisation $\lambda = \lambda'\lambda''
= \lambda''\gamma$ and $\lambda', \gamma \in {\mathcal
L}^{\le(M-k)}(E)$ (as $M-k \le \ell$).
We claim that
\[
Q = \sum_{\nu \in G} S_{\nu \lambda} P_{r ( [v]_\ell , \lambda )}
S_{\nu \lambda}^*
\]

\noindent is such that
\begin{align}
\| Q \pi_{S,P} ( \Phi (a) ) Q \| &= \| \pi_{S,P} ( \Phi (a) ) \| ,
\text{ and } \label{top} \\
Q S_\alpha P_{[v]_\ell} S_\beta^* Q &= 0 \text{ when } ( \alpha ,
[v]_\ell, \beta ) \in F \text{ and } \vert \alpha \vert \neq \vert
\beta \vert . \label{bottom}
\end{align}

\noindent The formula for $Q$ can be made  sense of by a calculation
similar to the one in Remark \ref{newspandef}. A routine calculation
verifies (\ref{top}).

Now suppose that $(\alpha , [v]_\ell , \beta) \in F$ satisfies
$\vert \alpha \vert \neq \vert \beta \vert$. Either $\alpha$ or
$\beta$ has length $k$, say $\vert \alpha \vert = k$.  As before,
$S_{\nu \lambda}^*S_\alpha$ is non-zero if and only if $\nu =
\alpha$.  Thus
\begin{eqnarray*}
Q S_\alpha P_{[v]_\ell} S_\beta^* Q &=& \sum_{ \nu \in G} S_{\alpha
\lambda}P_{r([v]_\ell,\lambda)} S_{\alpha \lambda}^* S_{\alpha }
P_{[v]_\ell}
S_{\beta}^*S_{\nu \lambda}P_{r([v]_\ell,\lambda)}S_{\nu \lambda}^*\\
 &=& \sum_{\nu
\in G} S_{\alpha \lambda}  P_{r([v]_\ell,\lambda )}
S_{\beta\lambda}^* S_{\nu \lambda} P_{r([v]_\ell,\lambda)} S_{\nu
\lambda}^* .
\end{eqnarray*}

\noindent For $P_{r([v]_\ell,\lambda)}S_{\beta\lambda}^* S_{\nu
\lambda}P_{r([v]_\ell,\lambda)}$ to be non-zero $\beta \lambda$ must
extend $\nu \lambda$, which implies that $\beta\lambda =
\nu\lambda\gamma$ for some $\gamma$. But then we have $\beta = \nu
\lambda'$ for some initial segment $\lambda'$ of $\lambda$ as $\vert
\beta \vert > \vert \nu \vert$. Hence  $\lambda = \lambda'
\lambda''$ which then implies that $\lambda=\lambda'' \gamma$ as
\[
\beta \lambda = \nu \lambda' \lambda'' = \nu \lambda \gamma = \nu
\lambda' \lambda'' \gamma
\]

\noindent and that $| \lambda' | = | \gamma |$. Since $| \beta | \le
M$ and $\vert \nu \vert = k$ it follows that $| \lambda' | \le M-k
\le \ell$. Thus $\lambda$ is agreeable for $[v]_\ell$, a
contradiction.
Thus $Q S_\alpha P_{[v]_\ell} S_\beta^* Q = 0$, and we have verified
(\ref{bottom}).

The rest of the proof is now standard (see, for example,
\cite[p.31]{cbms}).
\end{proof}

\section{Simplicity and Pure Infiniteness} \label{biggies}

Recall from \cite[Corollary 6.8]{kprr} that a directed graph $E$ is {\em cofinal} if for all
$x \in E^\infty$ and $v \in E^0$ there is a path $\lambda \in E^*$ and $N \ge 1$ such
that $s ( \lambda ) = v$ and $r ( \lambda ) = r ( x_N )$. Along with condition (L), cofinality is the key
hypothesis in the simplicity results for directed graphs (see \cite[Corollary 6.8]{kprr},
\cite[Proposition 5.1]{bprsz}).

In this section we seek an analogue for cofinality in the context of
labelled graphs which will allow us to prove a simplicity theorem
for labelled graph $C^*$-algebras. The two key difficulties to
overcome in the context of labelled graphs are that we must
accommodate the generalised vertices $[v]_\ell$ in a labelled graph
and the fact that there may be many representatives of a given infinite labelled path $x
\in \textsf{X}^+_{E,\mathcal{L}}$.

\begin{dfns}
Let $(E,\mathcal{L},\mathcal{E}^{0,-})$ be a labelled space and
$\ell \ge 1$. We say that $(E,\mathcal{L},\mathcal{E}^{0,-})$ is
{\em $\ell$-cofinal} if for all $x \in
\textsf{X}_{E,\mathcal{L}}^+$, $[v]_\ell \in \Omega_\ell$, and  $w
\in s(x)$ there is an $R(w)  \ge \ell$, an $N \ge 1$ and a finite
number of labelled paths $\lambda_1,\dots,\lambda_m$  such that for
all $r \ge R(w)$ we have $\bigcup_{i=1}^m r([v]_\ell,\lambda_i)
\supseteq r( [w]_r , x_1 \dots x_N )$.

We say that $(E,\mathcal{L},\mathcal{E}^{0,-})$ is {\em cofinal} if
there is an $L > 0$ such that $(E,\mathcal{L},\mathcal{E}^{0,-})$ is
{\em $\ell$-cofinal} for all $\ell > L$.
\end{dfns}

\begin{examples}
\begin{itemize}
\item[(i)] Recall from Example \ref{trivialex} that a directed graph $E$ may be
considered to be a labelled graph with the trivial labelling
$\mathcal{L}_t$. Let $E$ be a cofinal directed graph and fix $v \in
E^0$, $x \in E^\infty$. Since $w=s(x)$ is the only vertex with
$r(w,x_1  \ldots  x_n ) \neq \emptyset$ for all $n$, we may put $R(w) =1$ and
invoke cofinality of $E$ to get the required $N$ and $\lambda$ so
that $( E, \mathcal{L}_t, \mathcal{E}^{0,-} )$ is cofinal with $L =
1$. Thus the definition of cofinality for labelled graphs reduces to
the usual definition of cofinality for directed graphs.

\item[(ii)] The labelled space $(E_2, {\mathcal L}_2,{\mathcal E}^{0,-}_2 )$ of
Example \ref{agreex}(i) is not $\ell$-cofinal  for $\ell \ge 2$, and
so not cofinal. To see this, observe that $[w]_\ell = \{w\}$ for
$\ell \ge 2$ and there is no labelled path joining $w$ to the
infinite path $(100)^\infty$.

\end{itemize}
\end{examples}

\noindent The following result will allow us to prove cofinality for
many interesting examples.

\begin{lem} \label{transcofinal}  Let $(E,\mathcal{L},\mathcal{E}^{0,-})$ be a labelled space.
If $E$ is row-finite, transitive and ${\mathcal E}^{0,-}$ contains
$\{v\}$ for all $v \in E^0$ then
$(E,\mathcal{L},\mathcal{E}^{0,-})$ is cofinal with $L = 1$.
\end{lem}

\begin{proof}
Let $w \in E^0$.  Since $\{w\} \in {\mathcal E}^{0,-}$
there must be an $R(w) \ge 1$ such that $[w]_r = \{w\}$
for all $r \ge R(w)$.

Let $\ell \ge 1$ and choose $[v]_\ell \in \Omega_\ell$. Let $w \in
E^0$, and choose $R(w)$ as in the first paragraph. Let $x \in
\textsf{X}^+_{E,{\mathcal L}}$ be such that $w \in s(x)$. Let $N \ge
1$.  Then as $E$ is row-finite there are only finitely many paths
$\mu_1,\dots,\mu_m$ in $E$ with $s(\mu_i) = w$ and ${\mathcal
L}(\mu_i) = x_1 \dots x_N$. By transitivity of $E$ there are paths
$\lambda_1,\dots,\lambda_m \in E^*$ with $s(\lambda_i) = v$ and
$r(\lambda_i) = r(\mu_i)$. Then
\[
\bigcup_{i=1}^m r([v]_\ell, {\mathcal L}(\lambda_i)) \supseteq
r( [w]_r , x_1\dots x_N )
\]

\noindent
as required. Thus $(E , \mathcal{L} , \mathcal{E}^{0,-} )$  is cofinal with $L = 1$.
\end{proof}

\begin{thm} \label{simple}
Let $(E , \mathcal{L} , \mathcal{E}^{0,-} )$ be cofinal and
disagreeable.  Then $C^*(E , \mathcal{L} , \mathcal{E}^{0,-} )$ is
simple.
\end{thm}

\begin{proof}
Since every ideal in a $C^*$-algebra is the kernel of  a
representation, it suffices to prove that every non-zero
representation $\pi_{S,P}$ of $C^*(E , \mathcal{L} ,
\mathcal{E}^{0,-} )$ is faithful.    Suppose $\pi_{S,P}$ is a
non-zero  representation of $C^*(E , \mathcal{L} , \mathcal{E}^{0,-}
)$.  If we have $P_{[v]_\ell} = 0$ for all $v \in E^0$ and $\ell \ge
1$ then $\pi_{S,P} = 0$.  Thus there is a $w \in E^0$ and an $r \ge 1$
with $P_{[w]_r} \ne 0$. Fix $[v]_\ell \in \Omega_\ell$.  We aim to
prove that $P_{[v]_\ell} \ne 0$.  Since $[w]_r$ is the disjoint
union of finitely many equivalence classes $[w_i]_k$ whenever $k \ge
r$, for each $k$ there is an $i$ such that $P_{[w_i]_k} \ne 0$.  So
without loss of generality, for a given $[v]_\ell \in \Omega_\ell$,
we may assume that $r \ge R(w)$.

Since $(E,{\mathcal L},{\mathcal E}^{0,-})$ is set-finite we apply
(\ref{relationforvell}) of Proposition \ref{newspandef} to obtain
\[
P_{[w]_r} = \sum_{x_1 \in L^1_{[w]_r}} S_{x_1}
P_{r([w]_r,x_1)} S_{x_1}^* .
\]

\noindent Since the left-hand side is nonzero it follows that
$S_{x_1} P_{r([w]_r,x_1)} S_{x_1}^* \neq 0$ for some $x_1 \in L^1_{[w]_r}$
which implies that $P_{r([w]_r,x_1)} \neq 0$. Arguing as in the proof of
Proposition \ref{newspandef} we have
\[
P_{r([w]_r,x_1)} = \sum_{x_2 \in L^1_{r([w]_r,x_1)}} S_{x_2}
P_{r(r([w]_r,x_1),x_2)} S_{x_2}^*
\]

\noindent and so we may deduce that there is an $x_2$
with $P_{r(r([w]_r,x_1),x_2)} = P_{r([w]_r,x_1 x_2)}\ne 0$.
Continuing in this way we produce $x = x_1  x_2  \ldots \in
\textsf{X}^+_{E,\mathcal{L}}$ such that $P_{r([w]_r, x_1 \ldots x_n
)} \neq 0$ for all $n \ge 1$.

Let $\ell  \ge 1$ and $[v]_\ell \in \Omega_\ell$. Since $r > R(w)$, by
cofinality, there are finitely many labelled paths
$\lambda_1,\dots,\lambda_m$ and an $N \ge 1$ such that
$\bigcup_{i=1}^m r([v]_\ell,\lambda_i) \supseteq r([w]_r,x_1 \dots
x_N)$. Since $P_{r([w]_r,x_1,\dots,x_N)} \ne 0$ we must have
$P_{r([v]_\ell,\lambda_i)} \ne 0$ for some $i \in \{1,\dots,m\}$.
Since $r([v]_\ell,\lambda_i) \subseteq r(\lambda_i)$ it then follows
that $P_{r(\lambda_i)} \ne 0$ and hence
$S_{\lambda_i} \ne 0$. Since $P_{[v]_\ell} = \sum_{\lambda \in
L^{|\lambda_i|}_{[v]_\ell}} S_\lambda P_{r([v]_\ell,\lambda)}
S_\lambda^*$ it then follows that $P_{[v]_\ell} \ne 0$ as required.

Thus all the projections $P_ {[v]_\ell}$ are non-zero and
Theorem \ref{ckunique} implies that $\pi_{S,P}$ is
faithful, completing our proof.
\end{proof}

\begin{examples}  \label{agreecofinal}
\begin{itemize}
\item[(i)] The labelled space $( E_1 , {\mathcal L}_1 , {\mathcal
E}^{0,-}_1 )$, shown to be agreeable in Examples \ref{agreex}(i) is
cofinal with $L = 1$. This follows by Lemma \ref{transcofinal}(i)
since $E_1$ is row-finite, transitive and $\{ v \} \in {\mathcal
E}_1^{0,-}$ for all $v \in E_1^0$. Hence $C^* ( E_1 , {\mathcal L}_1
, \mathcal{E}^{0,-}_1 )$ is simple by Theorem \ref{simple}.

\item[(ii)] The labelled space $( E_{G,S} , \mathcal{L}_{G,S}, \mathcal{E}_{G,S}^{0,-})$
of Examples \ref{agreex}(ii) is cofinal with $L=1$.
To see this recall that $[g]_\ell = E_{G,S}^0=G$ for all $\ell \ge
1$. Fix $[g]_\ell \in \Omega_\ell$ and $x \in
\textsf{X}^+_{E_{G,S},\mathcal{L}_{G,S}}$. For $h \in G$, $r \ge
R(h)=1$ and $n=1$ we have $r ( [h]_r , x_1 ) = G$. Let $\lambda_1$
be any element of $S$, then $ r ( [g]_\ell , \lambda_1 ) = G  = r (
[h]_r , x_1 )$.  Hence $C^*(E_{G,S} , \mathcal{L}_{G,S},
\mathcal{E}_{G,S}^{0,-})$ is simple by Theorem \ref{simple}.
\end{itemize}
\end{examples}

\noindent We now turn our attention to the question of pure infiniteness for
simple labelled graph $C^*$-algebras. For graph $C^*$-algebras the key hypotheses
are condition (L) and every vertex connects to a loop (see \cite[Theorem 3.9]{kpr},
\cite[Proposition 5.4]{bprsz}). As we already have an
analogue of condition (L), we must now seek to find a suitable replacement for
the requirement that every vertex connects to a loop in the context of labelled graphs.
Again, there are two difficulties to overcome: we must accommodate the generalised
vertices $[v]_\ell$ in a labelled graph and find the correct analogue of a loop.

\begin{dfns}
The labelled path $\alpha$ is {\em repeatable} if $\alpha^n \in
{\mathcal L}^*(E)$ for all $n \ge 1$. We say that {\em every vertex
connects to a repeatable labelled path} if for every $[v]_m \in
\Omega_m$ there is a $w\in E^0$, $L(w) \ge 1$ and labelled paths
$\alpha,\delta \in {\mathcal L}^*(E)$ with $w \in
r([v]_m,\delta\alpha)$ such that $[w]_\ell \subseteq
r([w]_\ell,\alpha)$ for all $\ell \ge L(w)$.
\end{dfns}

\begin{rmk} \label{shouldbedef} The requirement that $[w]_\ell \subseteq r([w]_\ell,\alpha)$ for all
$\ell \ge L(w)$ ensures that $\alpha$ is repeatable, $\delta
\alpha^i \in \mathcal{L}^* ( E )$ for all $i \ge 1$ and that $r(
[w]_\ell, \alpha^i ) \ne \emptyset$ for all sufficiently large
$\ell$.
\end{rmk}

\noindent Our proof of the pure infiniteness result requires the following lemma
whose proof follows along similar lines to that of \cite[Lemma 5.4]{bprsz}.

\begin{lem} \label{getaproj}
Let $(E,{\mathcal L},{\mathcal E}^{0,-})$ be a labelled space, $v
\in E^0$ and $\ell \ge 1$.  Let $t$ be a positive element of
${\mathcal F}^k([v]_\ell)$.  Then there is a projection $r$ in the
$C^*$-subalgebra of ${\mathcal F}^k([v]_\ell)$ generated by $t$ such
that $rtr= \|t\|r$.
\end{lem}

\begin{thm} \label{pureinf}
Let $(E , \mathcal{L} , \mathcal{E}^{0,-} )$ be cofinal and
disagreeable.   If every vertex connects to a repeatable labelled
path then $C^*(E , \mathcal{L} , \mathcal{E}^{0,-} )$ is simple and
purely infinite.
\end{thm}

\begin{proof}
We know that $C^*(E,{\mathcal L},{\mathcal E}^{0,-})$ is simple by
Theorem \ref{simple}.  We show that every hereditary subalgebra
$A$ of $C^*(E,{\mathcal L},{\mathcal E}^{0,-})$ contains an infinite
projection; indeed we shall produce one which is dominated by a fixed
positive element $a \in A$ with  $\|\Phi(a)\| = 1$.

By Proposition \ref{newspandef} we may choose a positive element $b
\in \text{span}\{s_\alpha p_{[v]_\ell} s_\beta^* \;:\; \alpha, \beta
\in {\mathcal L}^*(E), [v]_\ell \in \Omega_\ell \}$  such that $\|a
- b\| < \frac{1}{4}$.  Suppose $b = \sum_{(\alpha , [w]_\ell,
\beta) \in F} c_{\alpha,[w]_\ell,\beta}
s_{\alpha}p_{[w]_\ell}s_{\beta}^*$ where $F$ is a finite subset of
$\mathcal{L}^* (E) \times \Omega \times \mathcal{L}^* (E)$.
The element $b_0 := \Phi(b)$  is positive and
satisfies $\|b_0\| \ge \frac{3}{4}$.

Let $k = \max \{ \vert \alpha \vert , \vert \beta \vert : ( \alpha ,
[w]_\ell , \beta ) \in F\}$. By repeatedly applying
\eqref{relationforvell} we may suppose (changing $F$ if necessary)
that $\min\{|\alpha|,|\beta| : ( \alpha , [w]_\ell , \beta ) \in F
\} = k$. Let $M = \max\{\vert \alpha \vert, \vert \beta \vert:
(\alpha , [w]_\ell , \beta) \in F\}$, $L_F = \max \{ L_w
\;:\; (\alpha,[w]_\ell,\beta) \in F\}$ and let $L$ be the smallest number
such that $(E,{\mathcal L},{\mathcal E}^{0,-})$ is $\ell$-cofinal for $\ell \ge L$.
Then from Proposition \ref{wasitfinitenoughforyou} and Remark \ref{wlog} we may assume
that $b_0 \in \oplus_{w : (\alpha,[w]_\ell,\beta) \in F }{\mathcal
F}^k([w]_m)$ for some $m \ge \max\{ L , L_F , M \}$. In fact, $\|b_0\|$ must be attained in some summand
${\mathcal F}^k([v]_m)$.  Let $b_1$ be the component of $b_0$ in
${\mathcal F}^k([v]_m)$, and note that $b_1 \ge 0$ and $\|b_1\| =
\|b_0\|$. By Lemma \ref{getaproj} there is a projection $r \in
C^*(b_1) \subseteq {\mathcal F}^k([v]_m)$ such that
$rb_1r=\|b_1\|r$. Since $b_1$ is a finite sum of
$s_{\alpha}p_{[v]_m}s_{\beta}^*$ we can write $r$ as a sum $\sum
c_{\alpha\beta}s_\alpha p_{[v]_m}s_\beta^*$ over all pairs of paths
in
\[
S = \{ \alpha \in {\mathcal L} ( E^k ) \; : \; \text{either } ( \alpha ,
[w]_\ell, \beta ) \in F \text{ or } ( \beta , [w]_\ell , \alpha )
\in F \text{ and } [w]_\ell \subseteq r(\alpha) \}.
\]

\noindent As $m \ge L_v$, $[v]_m$ is disagreeable and there is an $n > M$ and
a $\lambda \in {\mathcal L}^*(E)$ with $|\lambda| \ge n$ which is
disagreeable for $[v]_m$.  Since $m \ge M \ge M-k$ as well we may
employ the same argument as in the proof of Theorem \ref{ckunique}
to produce a projection $Q := \sum_{\gamma \in S} s_{\gamma\lambda}
p_{r([v]_m,\lambda)}s_{\gamma\lambda}^*$ such that $Q
s_{\alpha}p_{[v]_m}s_{\beta}^*Q = 0$ unless $|\alpha| = |\beta| = k$
and $[v]_m \subseteq r(\alpha) \cap r(\beta)$. Since $r \in
C^*(b_1)$ we have
\[
r= \sum c_{\alpha\beta} s_\alpha p_{[v]_m} s_\beta^* = \sum
c_{\alpha\beta}s_\alpha(s_{\lambda}p_{r([v]_m,\lambda)}s_\lambda^* +
(p_{[v]_m} - s_\lambda
 p_{r([v]_m,\lambda)}s_\lambda^*))s_\beta^* \ge Q
\]

\noindent so that
\[
QbQ = Qb_0Q = Qrb_1rQ = \|b_1\|rQ = \|b_0\|Q \ge \frac{3}{4}Q.
\]

\noindent Since $\|a - b\| \le \frac{1}{4}$ we have $QaQ \ge QbQ -
\frac{1}{4}Q \ge \frac{1}{2}Q$ and so $QaQ$ is invertible in
$QC^*(E,{\mathcal L},{\mathcal E}^{0,-})Q$.  Let $c$ denote its
inverse and put $v = c^{1/2}Qa^{1/2}$.  Then  $vv^* =
c^{1/2}QaQc^{1/2} = Q$, and $v^*v = a^{1/2}QcQa^{1/2} \le \|c\|a$
and so $v^*v$ belongs to the hereditary subalgebra $A$. To finish,
we must show that $v^*v$ is an infinite projection.

We wish to find a labelled path $\beta$ with $r([v]_m,\beta) \ne \emptyset$
whose initial segment is $\lambda$ and whose terminal segment
is a repeatable labelled path.   We choose $x \in
r([v]_m,\lambda)$.  Then $[x]_{m + |\lambda|} \subseteq
r([v]_m,\lambda)$ and by hypothesis $[x]_{m+\vert \lambda \vert}$ connects to a repeatable
path: That is, there is a $w \in E^0$, $L(w) \ge 1$ and paths
$\alpha, \delta \in {\mathcal L}^*(E)$ such that $w \in r([x]_{m +
|\lambda|},\delta \alpha)$,  and $[w]_n \subseteq r([w]_n,\alpha)$
for all $n \ge L(w)$.   The required path is $\beta = \lambda \delta
\alpha$.   Let $N = \max\{L_w,L(w)\}$.  We claim that
$p_{[w]_n}$ is an infinite projection for all $n \ge N$. As $n \ge
L(w)$, we know from Remark \ref{shouldbedef} that we have $r ( [ w
]_n , \alpha^i ) \neq \emptyset$, for $i \ge 1$. Moreover, as $n \ge
L_w$ we know that $[w]_n$ is disagreeable. Hence there must be a
labelled path $\gamma$ with $[w]_n \cap s ( \gamma ) \neq \emptyset$
and $i \ge 1$ with $|\gamma| = |\alpha^i |$, and $\gamma \ne
\alpha^i$. We compute
\[
p_{[w]_n} \le s_{\alpha^i} p_{r([w]_n,\alpha^i)} s_{\alpha^i}^*
              < s_{\alpha^i} p_{r([w]_n,\alpha^i)}s_{\alpha^i}^* +
                s_\gamma p_{r([w]_n,\gamma)}s_\gamma^* \le
p_{[w]_n}
\]

\noindent and our claim is established.

We now demonstrate the existence of an infinite subprojection of
$Q$. If $\mu$ is any labelled path with $|\mu| = k \le M \le m$ and $r(\mu)
\cap s(\lambda) \cap [v]_m \ne \emptyset$ then for $n \ge N$ such
that $[w]_n \subseteq r([v]_m,\lambda\delta\alpha)$ (note that such
an $n$ exists as $[w]_n \subseteq r([v]_m,\lambda\delta \alpha)$
for all sufficiently large $n$) we have
\[
p_{[w]_n} = p_{[w]_n} s_{\mu\lambda \delta \alpha}^* s_{\mu\lambda
\delta \alpha} \thicksim s_{\mu\lambda\delta\alpha}p_{[w]_n}
s_{\mu\lambda\delta\alpha}^* \le s_{\mu\lambda}p_{r([v]_m,\lambda)}
s_{\mu\lambda}^* .
\]

\noindent Because the projection $s_{\mu\lambda}p_{r([v]_m,\lambda)}
s_{\mu\lambda}^*$ is a minimal projection in the matrix algebra
$\text{span}\{s_{\mu\lambda}p_{r([v]_m,\lambda)} s_{\nu\lambda}^*
\;:\; \mu,\nu \in S\}$, it is equivalent to a subprojection of $Q$.
It follows that $Q$ is infinite, and, since $Q = vv^* \thicksim
v^*v$ this completes the proof.
\end{proof}

\begin{examples}
\begin{itemize}
\item[(i)]  In the labelled space $(E_1,{\mathcal L}_1,{\mathcal
E}^{0,-}_1)$ of Examples \ref{agreex}(i) every vertex connects to
the repeatable path $0$.  Since $(E_1,{\mathcal L}_1, {\mathcal
E}^{0,-}_1)$ is cofinal and disagreeable, $C^*(E_1,{\mathcal
L}_1,{\mathcal E}^{0,-}_1)$ is simple and purely infinite by Theorem
\ref{pureinf}.

\item[(ii)] Suppose that for a group $G$, the set $S$ contains
(not necessarily distinct) elements $g_1,\dots,g_n$ such that $g_1
\dots g_n = 1_G$. Then every vertex in the labelled graph $( E_{G,S}
, \mathcal{L}_{G,S})$ of Examples \ref{agreex}(ii) connects to the
repeatable labelled path $g_1 \dots g_n$.  If in addition we have
$|S|>1$, then by Examples \ref{agreex} (ii) and Examples
\ref{agreecofinal} (ii)  $( E_{G,S},\mathcal{L}_{G,S}, {\mathcal
E}_{G,S}^{0,-})$ is cofinal and disagreeable and so
$C^*(E_{G,S},\mathcal{L}_{G,S}, {\mathcal E}_{G,S}^{0,-})$ is simple
and purely infinite by Theorem \ref{pureinf}.
\end{itemize}
\end{examples}

\section{Some labelled graph presentations of non-sofic shift spaces} \label{newex}

\subsection{Dyck Shifts} \label{nogashere}

In this section we associate a labelled
graph to a Dyck shift in such a way that the resulting labelled
space $C^*$-algebra is simple and purely infinite.

First we recall the definition of the Dyck shift (see, for example,
\cite{mey1,ma2}). Let $N \ge 1$ be a fixed positive integer. The
Dyck shift $D_N$ has alphabet ${\mathcal A} =
\{\alpha_1,\dots,\alpha_N,\beta_1,\dots,\beta_N\}$ where the symbols
$\alpha_i$ correspond to opening brackets of type $i$ and
the symbols $\beta_i$ are their respective closing brackets.
We say that a word $\gamma_1\dots\gamma_n \in {\mathcal A}^*$ is
admissible if  $\gamma_1\dots\gamma_n$ does not contain any
substring $\alpha_i\beta_j$ with $i \neq j$.  Thus the language of
the Dyck shift consists of all strings of properly matched brackets
of types $\alpha_1,\dots \alpha_N$.

The following algorithm gives a
labelled graph presentation of a Dyck shift.

\begin{enumerate}
\item Fix $N \ge 1$ and an alphabet $\{\alpha_1,\dots,\alpha_N,\beta_1,\dots,\beta_N\}$.
\item Draw an unrooted, infinite, directed  tree in which every vertex receives
one edge and emits $N$ edges (i.e. an $N$-ary tree).  Label the $N$ branches from each
node, working from left to right, by $\alpha_1,\dots,\alpha_N$.
\item  For each $i \in \{1,\dots,N\}$ and each edge $e$ labelled
$\alpha_i$, draw an edge from $r(e)$ to $s(e)$ with label $\beta_i$.
\end{enumerate}

\noindent
The resulting labelled graph $(E_N,{\mathcal L}_N)$ is a
left-resolving labelled graph which presents the
Dyck shift $D_N$.

\begin{examples}
\begin{enumerate}
\item Let $N =1$ and ${\mathcal A} = \{ \;(\;,\;)\;\}$. The above algorithm gives the following labelled
graph presentation of $D_1$.
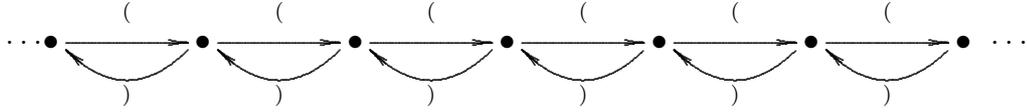
\[
\beginpicture
\setcoordinatesystem units <2cm,1cm>

\setplotarea x from -3.1 to 3, y from -0.3 to 0.3

\put{$\ldots$}[l] at -3.3 0

\put{$\bullet$} at -2 0

\put{$\bullet$} at -1 0

\put{$\bullet$} at 0 0

\put{$\bullet$} at 1 0

\put{$\bullet$} at 2 0

\put{$\bullet$} at -3 0

\put{$\bullet$} at 3 0

\put{$\ldots$} at 3.3 0

\put{{\tiny $($}} at -2.5 0.4

\put{{\tiny $($}} at -1.5 0.4

\put{{\tiny $($}} at -0.5 0.4

\put{{\tiny $($}} at 0.5 0.4

\put{{\tiny $($}} at 1.5 0.4

\put{{\tiny $($}} at 2.5 0.4

\put{{\tiny $)$}} at -2.5 -0.7

\put{{\tiny $)$}} at -1.5 -0.7

\put{{\tiny $)$}} at -0.5 -0.7

\put{{\tiny $)$}} at 0.5 -0.7

\put{{\tiny $)$}} at 1.5 -0.7

\put{{\tiny $)$}} at 2.5 -0.7

\setquadratic

\plot -1.1 -0.1 -1.5 -0.5 -1.9 -0.1 /

\plot -2.1 -0.1 -2.5 -0.5 -2.9 -0.1 /

\plot -0.1 -0.1 -0.5 -0.5 -0.9 -0.1 /

\plot 0.1 -0.1 0.5 -0.5 0.9 -0.1 /

\plot 1.1 -0.1 1.5 -0.5 1.9 -0.1 /

\plot 2.1 -0.1 2.5 -0.5 2.9 -0.1 /

\arrow <0.25cm> [0.1,0.3] from -2.9 0 to -2.1 0

\arrow <0.25cm> [0.1,0.3] from -1.9 0 to -1.1 0

\arrow <0.25cm> [0.1,0.3] from -0.9 0 to -0.1 0

\arrow <0.25cm> [0.1,0.3] from 0.1 0 to 0.9 0

\arrow <0.25cm> [0.1,0.3] from 1.1 0 to 1.9 0

\arrow <0.25cm> [0.1,0.3] from 2.1 0 to 2.9 0

\arrow <0.25cm> [0.1,0.3] from -2.8 -0.27 to -2.9 -0.1

\arrow <0.25cm> [0.1,0.3] from -1.8 -0.27 to -1.9 -0.1

\arrow <0.25cm> [0.1,0.3] from -0.8 -0.27 to -0.9 -0.1

\arrow <0.25cm> [0.1,0.3] from 0.2 -0.27 to 0.1 -0.1

\arrow <0.25cm> [0.1,0.3] from 1.2 -0.27 to 1.1 -0.1

\arrow <0.25cm> [0.1,0.3] from 2.2 -0.27 to 2.1 -0.1
\endpicture
\]

\noindent Of course, the above labelled graph is not the optimal presentation of $D_1$, as $D_1$ 
has no constraints and so is the full shift on the symbols $($ and $)$.

\item Let $N = 2$ and let ${\mathcal A} = \{\;(,\;[,\;),\;]\;\}$.  The above algorithm gives
the following labelled graph presentation of $D_2$.
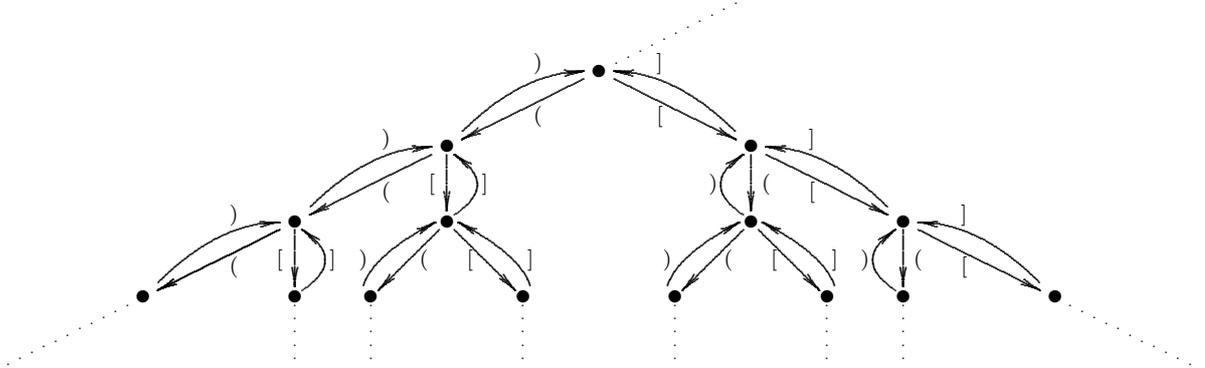
\[
\beginpicture

\setcoordinatesystem units <2cm,1cm>

\setplotarea x from -3.1 to 3, y from -2.1 to 3

\put{$\bullet$} at  0 2

\put{$\bullet$} at -1 1

\put{$\bullet$} at 1  1

\put{$\bullet$} at -1 0

\put{$\bullet$} at -2 0

\put{$\bullet$} at 1 0

\put{$\bullet$} at 2 0

\put{$\bullet$} at -1.5 -1

\put{$\bullet$} at -2 -1

\put{$\bullet$} at -0.5 -1

\put{$\bullet$} at 0.5 -1

\put{$\bullet$} at 1.5 -1

\put{$\bullet$} at 2 -1

\put{$\bullet$} at -3 -1

\put{$\bullet$} at 3 -1

\put{{\tiny $($}} at -0.4 1.4

\put{{\tiny $($}} at -1.4 0.4

\put{{\tiny $($}} at -2.4 -0.6

\put{{\tiny $[$}} at 0.4 1.4

\put{{\tiny $[$}} at 1.4 0.4

\put{{\tiny $[$}} at 2.4 -0.6

\put{{\tiny $)$}} at -0.4 2.1

\put{{\tiny $)$}} at -1.4 1.1

\put{{\tiny $)$}} at -2.4 0.1

\put{{\tiny $]$}} at 0.4 2.1

\put{{\tiny $]$}} at 1.4 1.1

\put{{\tiny $]$}} at 2.4 0.1

\put{{\tiny $[$}} at -1.1 0.5

\put{{\tiny $]$}} at -0.75 0.5

\put{{\tiny $[$}} at -2.1 -0.5

\put{{\tiny $]$}} at -1.75 -0.5

\put{{\tiny $($}} at 1.1 0.5

\put{{\tiny $)$}} at 0.75 0.5

\put{{\tiny $($}} at 2.1 -0.5

\put{{\tiny $)$}} at 1.75 -0.5

\put{{\tiny $($}} at -1.15 -0.5

\put{{\tiny $)$}} at -1.55 -0.5

\put{{\tiny $($}} at 0.85 -0.5

\put{{\tiny $)$}} at 0.45 -0.5

\put{{\tiny $]$}} at -0.45 -0.5 

\put{{\tiny $[$}} at -0.85 -0.5

\put{{\tiny $]$}} at 1.55 -0.5

\put{{\tiny $[$}} at 1.15 -0.5

\arrow <0.25cm> [0.1,0.3] from -0.1 1.9 to -0.9 1.1

\arrow <0.25cm> [0.1,0.3] from 0.1 1.9 to 0.9 1.1

\arrow <0.25cm> [0.1,0.3] from -1.1 0.9 to -1.9 0.1

\arrow <0.25cm> [0.1,0.3] from -2.1 -0.1 to -2.9 -0.9

\arrow <0.25cm> [0.1,0.3] from -2.1 -0.1 to -2.9 -0.9

\arrow <0.25cm> [0.1,0.3] from 1.1 0.9 to 1.9 0.1

\arrow <0.25cm> [0.1,0.3] from 2.1 -0.1 to 2.9 -0.9

\arrow <0.25cm> [0.1,0.3] from -1 0.9 to -1 0.15

\arrow <0.25cm> [0.1,0.3] from 1 0.9 to 1 0.15

\arrow <0.25cm> [0.1,0.3] from -2 -0.1 to -2 -0.85

\arrow <0.25cm> [0.1,0.3] from 2 -0.1 to 2 -0.85

\arrow <0.25cm> [0.1,0.3] from 1.05 -0.1 to 1.45 -0.9

\arrow <0.25cm> [0.1,0.3] from -0.95 -0.1 to -0.55 -0.9

\arrow <0.25cm> [0.1,0.3] from -1.05 -0.1 to -1.45 -0.9

\arrow <0.25cm> [0.1,0.3] from 0.95 -0.1 to 0.55 -0.9

\setquadratic

\plot -0.1 2   -0.5  1.8    -0.9 1.2 /

\plot -1.1 1   -1.5  0.8    -1.9 0.2 /

\plot -2.1 0   -2.5  -0.2    -2.9 -0.8 /

\plot 0.1 2   0.5  1.8    0.9 1.2 /

\plot 1.1 1   1.5  0.8    1.9 0.2 /

\plot 2.1 0   2.5  -0.2    2.9 -0.8 /

\plot  -0.95 0.9 -0.8 0.5  -0.95 0.1 /

\plot  0.95 0.9 0.8 0.5  0.95 0.1 /

\plot  -1.95 -0.1 -1.8 -0.5 -1.95 -0.9 /

\plot  1.95 -0.1 1.8 -0.5 1.95 -0.9 /

\plot 1.075 0 1.4 -0.4275 1.55 -0.85 /

\plot -0.925 0 -0.6 -0.4275 -0.45 -0.85 /

\plot -1.075 0 -1.4 -0.4275 -1.55 -0.85 /

\plot 0.925 0 0.6 -0.4275 0.45 -0.85 /

\arrow <0.25cm> [0.1,0.3] from -0.2 1.98 to -0.1 2

\arrow <0.25cm> [0.1,0.3] from -1.2 0.98 to -1.1 1

\arrow <0.25cm> [0.1,0.3] from -2.2 -0.02 to -2.1 0

\arrow <0.25cm> [0.1,0.3] from 0.2 1.98 to 0.1 2

\arrow <0.25cm> [0.1,0.3] from 1.2 0.98 to 1.1 1

\arrow <0.25cm> [0.1,0.3] from 2.2 -0.02 to 2.1 0

\arrow <0.25cm> [0.1,0.3] from -0.885 0.8 to -0.95 0.9

\arrow <0.25cm> [0.1,0.3] from -1.885 -0.2 to -1.95 -0.1

\arrow <0.25cm> [0.1,0.3] from 0.885 0.8 to 0.95 0.9

\arrow <0.25cm> [0.1,0.3] from 1.885 -0.2 to 1.95 -0.1

\arrow <0.25cm> [0.1,0.3] from 1.165 -0.1 to 1.075 0

\arrow <0.25cm> [0.1,0.3] from -0.835 -0.1 to -0.925 0

\arrow <0.25cm> [0.1,0.3] from 0.835 -0.1 to 0.925 0

\arrow <0.25cm> [0.1,0.3] from -1.165 -0.1 to -1.075 0

\setdots

\arrow <0cm> [0,0] from 0.9 2.9 to  0.1 2.1

\arrow <0cm> [0,0] from -3.1 -1.1 to  -3.9 -1.9

\arrow <0cm> [0,0] from -2 -1.1 to  -2 -1.9

\arrow <0cm> [0,0] from -1.5 -1.1 to -1.5 -1.9

\arrow <0cm> [0,0] from -0.5 -1.1 to -0.5 -1.9

\arrow <0cm> [0,0] from 0.5 -1.1 to 0.5 -1.9

\arrow <0cm> [0,0] from 1.5 -1.1 to 1.5 -1.9

\arrow <0cm> [0,0] from 2 -1.1 to 2 -1.9

\arrow <0cm> [0,0] from 3.1 -1.1 to 3.9 -1.9

\endpicture
\]
\end{enumerate}
\end{examples}

\begin{prop} \label{dyckprop}
Let $N \ge 1$ and  ${\mathcal
A}=\{\alpha_1,\dots,\alpha_N,\beta_1,\dots,\beta_N\}$. Then
$C^*(E_N,{\mathcal L}_N,{\mathcal E}_N^{0,-})$ is simple and purely
infinite.
\end{prop}

\begin{proof}
Let $\Lambda_\ell = \{ \lambda^\ell_1 , \ldots ,
\lambda^\ell_{N^\ell} \}$ be the labelled paths of length $\ell$
which consist of only $\alpha_i$'s (opening braces), and let
$\Xi_\ell = \{ \mu^\ell_1 , \ldots , \mu^\ell_{N^\ell} \}$ be the
labelled paths of length $\ell$ which consist of only $\beta_i$'s
(closing braces), organised in such a way that  for all $i$ the word
$\lambda^\ell_i \mu^\ell_i$ belongs to the language of $D_N$. Since
every vertex $v \in E_N^0$ receives one opening brace and $N$
closing braces, it follows that $v$ receives a unique
$\lambda^\ell_i \in \Lambda_\ell$ one sees that $\Omega_\ell = \{ [
v^\ell_i ]_\ell : i = 1 , \ldots , N^\ell \}$ where $v^\ell_i$ is
some vertex in $r ( \lambda^\ell_i )$.  Moreover, every vertex $v
\in E_N^0$ emits exactly one closing brace (the closing version of
the one it receives) and $N$ opening braces, so every $v$ which
receives $\lambda^\ell_i$ also emits $\mu^\ell_i$.

For $1 \le i , j \le N^\ell$ let $\mu^\ell_{ij} = \mu^\ell_i \lambda^\ell_j$ then
$s ( \mu^\ell_{ij} ) = [ v^\ell_i ]_\ell$ as the only vertices which emit $\mu^\ell_i$ are those which
receive $\lambda^\ell_i$. Moreover, we have $r ( \mu^\ell_{ij} ) = r ( \lambda^\ell_j ) = [ v^\ell_j ]_\ell$ since
every vertex in $E_N^0$ (emits the labelled path $\lambda^\ell_i$ and hence) receives a
labelled path $\mu^\ell_i$ which originates from a vertex in $[ v^\ell_i ]_\ell$,
that it $r ( [v^\ell_i ]_\ell , \mu^\ell_i ) = E_N^0$.

Fix $\ell \ge 1$, $[v]_\ell \in \Omega_\ell$ and $x \in \textsf{X}_{E_N,\mathcal{L}_N}^+$.
Without loss of generality suppose that $[v]_\ell = [v^\ell_1 ]_\ell$.
Then by definition of the $\mu^\ell_{ij}$ we have
\[
\bigcup_{j=1}^{N^\ell} r ( [v^\ell_1 ]_\ell , \mu^\ell_{1j} ) =
E_N^0
\]

\noindent and hence the labelled space
$(E_N,{\mathcal L}_N, {\mathcal E}_N^{0,-})$ is cofinal with $L=1$.

We now show that $(E_N , {\mathcal L}_N, {\mathcal E}_N^{0,-})$ is disagreeable.
For $n \ge 1$, every vertex $v$ emits the labelled path $\alpha_1^n \beta_1$, which
is disagreeable for $[v]_\ell$. Hence $[v]_\ell$ is disagreeable for all $\ell \ge 1$.
It follows that $C^*(E_N , {\mathcal L}_N , {\mathcal E}_N^{0,-})$
is simple by Theorem \ref{simple}.

Since every vertex $v \in E_N^0$ emits the repeatable labelled path
$\alpha_1 \beta_1$  it follows that every generalised vertex in $(
E_N , {\mathcal L}_N, {\mathcal E}_N^{0,-} )$ connects to a
repeatable labelled path.  Thus $C^*( E_N , {\mathcal L}_N ,
{\mathcal E}_N^{0,-})$ is purely infinite by Theorem \ref{pureinf}.
\end{proof}

\begin{rmk} \label{matsumotoagain}
The essential symbolic matrix system $( M(E_N)_{\ell,\ell+1} ,
I(E_N)_{\ell,\ell+1}  )_{\ell \ge 1}$ associated to $( E_N , {\mathcal
L}_N , {\mathcal E}_N^{0,-})$ gives rise to the $\lambda$-graph system
$\mathfrak{L}^{Ch(D_N)}$ described on \cite[p.5]{ma2} (for example).
Hence by Proposition \ref{wegotmat} it follows that $C^*
( E_N , {\mathcal L}_N , {\mathcal E}_N^{0,-} ) \cong
\mathcal{O}_{\mathfrak{L}^{Ch(D_N)}}$. Moreover by \cite[Proposition
5.1]{ma2} we know that
\[
K_0 ( C^* ( E_N , {\mathcal L}_N , {\mathcal E}_N^{0,-}) ) \cong
{\bf Z} / N {\bf Z} \oplus C ( \mathfrak{K} , {\bf Z} )
\text{ and } K_1 ( C^* ( E_N , {\mathcal L}_N , {\mathcal E}_N^{0,-}) ) \cong 0
\]

\noindent where $C ( \mathfrak{K} , {\bf Z} )$ denotes the abelian
group of all ${\bf Z}$-valued continuous functions on the Cantor set
$\mathfrak{K}$. Since the $K$-theory of $C^* ( E_N , {\mathcal
L}_N , {\mathcal E}_N^{0,-} )$ is not finitely generated it follows that $C^*
( E_N , {\mathcal L}_N , {\mathcal E}_N^{0,-} )$ cannot be isomorphic to a
unital graph algebra (indeed $C^* ( E_N , {\mathcal L}_N , {\mathcal
E}_N^{0,-} )$ is not semiprojective).

Note that the essential symbolic matrix system $(
M (E_N)_{\ell,\ell+1} , I (E_N)_{\ell,\ell+1}  )_{\ell \ge 1}$
associated to $( E_N , {\mathcal L}_N , {\mathcal E}_N^{0,-} )$ is not the
same as the one described in \cite[Proposition 2.1]{km}. In
\cite{km} the $\lambda$-graphs associated to symbolic matrix systems
are ``upward directed" whereas in \cite{m99} they are ``downward
directed". This results from the change of time direction mentioned
on \cite[p.81]{km}. Hence to form the appropriate ``upward directed"
versions for $( E_N , {\mathcal L}_N , {\mathcal E}_N^{0,-})$, it would
seem natural to reverse the arrows in $E_N$.
\end{rmk}

\subsection{A Further Example} \label{wannabenogashere}

Consider the shift space $X$ over the alphabet ${\mathcal A} = \{a,b,c\}$
whose language consists of all words in $\{a,b,c\}$ such that the numbers of
$b$'s and $c$'s occurring between any pair of consecutive $a$'s are equal.

Note that the  shift $X$ is not sofic:
suppose otherwise. Then there is a finite labelled graph
$(E_X,{\mathcal L}_X)$ with $|E_X^0| = n$ which presents $X$. Let $\alpha$ be a
path in $(E_X,{\mathcal L}_X)$ which presents
$ab^{2n}c^{2n}a$. Then since the number of $c$'s in ${\mathcal
L}_X(\alpha)$ is greater than $n$, $\alpha$ must contain
a cycle $\tau$ such that ${\mathcal L}_X(\tau) = c^m$ for
some $m \le n$. Write $\alpha = \alpha'\tau\alpha''$. Then $\beta =
\alpha' \tau^2 \alpha''$ is a path in $E_X^*$ which
presents the forbidden word $ab^{2n}c^{2n+m}a$.

The  shift $X$ has the following labelled
graph presentation $(E_X,{\mathcal L}_X)$:
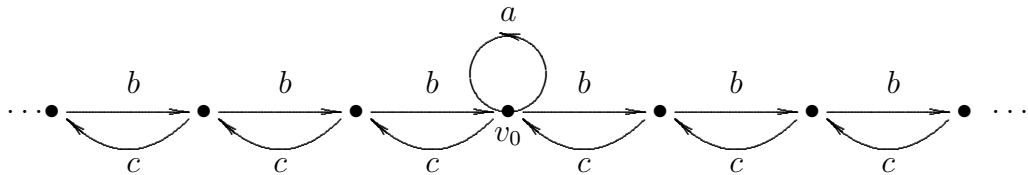
\[
\beginpicture
\setcoordinatesystem units <2cm,1cm>

\setplotarea x from -3.1 to 3, y from  0.3  to 2

\put{$\ldots$}[l] at -3.3 0

\put{$\bullet$} at -2 0

\put{$\bullet$} at -1 0

\put{$\bullet$} at 0 0

\put{$\bullet$} at 1 0

\put{$\bullet$} at 2 0

\put{$\bullet$} at -3 0

\put{$\bullet$} at 3 0

\put{$\ldots$} at 3.3 0

\put{ $b$} at -2.5 0.4

\put{ $b$} at -1.5 0.4

\put{$b$} at -0.5 0.4

\put{$b$} at 0.5 0.4

\put{$b$} at 1.5 0.4

\put{ $b$} at 2.5 0.4

\put{ $c$} at -2.5 -0.7

\put{$c$} at -1.5 -0.7

\put{$c$} at -0.5 -0.7

\put{$c$} at 0.5 -0.7

\put{$c$} at 1.5 -0.7

\put{$c$} at 2.5 -0.7

\put{$a$}[b] at 0 1.2

\put{$v_0$}[t] at 0 -0.2

\circulararc 360 degrees from 0 0 center at 0 0.5

\setquadratic

\plot -1.1 -0.1 -1.5 -0.5 -1.9 -0.1 /

\plot -2.1 -0.1 -2.5 -0.5 -2.9 -0.1 /

\plot -0.1 -0.1 -0.5 -0.5 -0.9 -0.1 /

\plot 0.1 -0.1 0.5 -0.5 0.9 -0.1 /

\plot 1.1 -0.1 1.5 -0.5 1.9 -0.1 /

\plot 2.1 -0.1 2.5 -0.5 2.9 -0.1 /

\arrow <0.25cm> [0.1,0.3] from 0.05 1 to -0.05 1.04

\arrow <0.25cm> [0.1,0.3] from -2.9 0 to -2.1 0

\arrow <0.25cm> [0.1,0.3] from -1.9 0 to -1.1 0

\arrow <0.25cm> [0.1,0.3] from -0.9 0 to -0.1 0

\arrow <0.25cm> [0.1,0.3] from 0.1 0 to 0.9 0

\arrow <0.25cm> [0.1,0.3] from 1.1 0 to 1.9 0

\arrow <0.25cm> [0.1,0.3] from 2.1 0 to 2.9 0

\arrow <0.25cm> [0.1,0.3] from -2.8 -0.27 to -2.9 -0.1

\arrow <0.25cm> [0.1,0.3] from -1.8 -0.27 to -1.9 -0.1

\arrow <0.25cm> [0.1,0.3] from -0.8 -0.27 to -0.9 -0.1

\arrow <0.25cm> [0.1,0.3] from 0.2 -0.27 to 0.1 -0.1

\arrow <0.25cm> [0.1,0.3] from 1.2 -0.27 to 1.1 -0.1

\arrow <0.25cm> [0.1,0.3] from 2.2 -0.27 to 2.1 -0.1

\endpicture
\]

\noindent Since the graph $E_X$ is transitive, it is straightforward to check from the
above presentation that $X$ is irreducible.

Since each vertex in $E_X$ to the right (resp.\ left) of $v_0$
receives a unique labelled path of the form $a b^n$ (resp.\ $ac^n$)
it follows that $\{ v \} \in {\mathcal E}^{0,-}$ for all $v \in
E_X^0$. Since $E_X$ is row-finite it follows that $( E_X ,
\mathcal{L}_X , \mathcal{E}_X^{0,-} )$ is cofinal by Lemma
\ref{transcofinal}.

For $n \ge 1$ every $v \in E_X^0$ emits the labelled path $b^nc$, which
is disagreeable for $[v]_\ell$. Hence $[v]_\ell$ is disagreeable for all $\ell \ge 1$ and
so $C^* ( E_X , \mathcal{L}_X , \mathcal{E}_X^{0,-} )$
is simple by Theorem \ref{simple}.

Every $v \in E_X^0$ emits the repeatable path $bc$ and since $E_X$ is transitive,
it follows that every generalised vertex
connects to a repeatable path.  Thus $C^* ( E_X , \mathcal{L}_X , \mathcal{E}_X^{0,-} )$
is purely infinite by Theorem
\ref{pureinf}.


\begin{thebibliography}{20}
\bibitem{bprsz} T.~ Bates, D.~ Pask, I.~ Raeburn and W.~Szyma\'{n}ski.
\newblock {\em The $C^*$-algebras of row-finite graphs.}
\newblock New York J.\ Math. {\bf 6} (2000), 307--324.

\bibitem{bp} T.\ Bates and D.\ Pask. {\em Flow equivalence of graph
algebras}.  Ergod.\ Th.\ \& Dynam.\ Sys., {\bf 24} (2004), 367--382.

\bibitem{bp2} T.\ Bates and D.\ Pask. {\em $C^*$-algebras of
labelled graphs}. J.\ Operator Theory {\bf 57} (2007), 207--226.

\bibitem{cm} T.\ Carlsen and K.\ Matsumoto. {\em Some remarks on
the $C^*$-algebras associated with subshifts}. Math.\ Scand.\ {\bf
95} (2004), 145--160.

\bibitem{ce} T.\ Carlsen and S.\ Eilers. {Matsumoto $K$-groups
associated to certain shift spaces}. Doc.\ Math.\ {\bf 9}
(2004),639--671.

\bibitem{cs} T.\ Carlsen and S.\ Silvestrov. {\em $C^*$-crossed products and shift spaces}.
preprint, arXiv: math.OA/0505503.

\bibitem{tc} T.\ Carlsen.  {\em Symbolic dynamics, partial dynamical systems, boolean
algebras and $C^*$-algebras generated by partial isometries}.
preprint, arXiv: math.OA/0604165.

\bibitem{km} W.~Krieger and K.~Matsumoto. \newblock{A Lambda-Graph system for the Dyck Shift and Its $K$-Groups}.
Doc.\ Math.\ {\bf 8} (2003), 79-96.

\bibitem{kprr} A.~Kumjian, D.~Pask, I.~Raeburn, and J.~Renault.
\newblock {\em Graphs, groupoids and Cuntz--Krieger algebras.}
\newblock J. Funct. Anal.  {\bf 144} (1997), 505--541.

\bibitem{kpr} A. Kumjian, D. Pask and I. Raeburn, {\em Cuntz-Krieger algebras of directed graphs}, Pacific J. Math. {\bf 184} (1998), 161--174.

\bibitem{kp}  A. \ Kumjian and D. \ Pask. {\em $C^*$-algebras of directed graphs and group actions}.  Ergod. \ Th. \ \& Dynam. \ Sys. \ {\bf 19} (1999), 1503--1519.

\bibitem{lm} D. Lind and B. Marcus. {\em An Introduction to Symbolic Dynamics and Coding}, CUP, 1995.

\bibitem{m} K.\ Matsumoto. {\em On $C^*$-algebras associated with
subshifts}. Internat.\ J.\ Math.\ {\bf 8}, (1997), 357-374.

\bibitem{ma3} K.\ Matsumoto. {\em Dimension groups for subshifts and simplicity
of the associated $C^*$-algebras}. J.\ Math.\ Soc.\ Japan {\bf 51}
(1999), 679--697.

\bibitem{m99} K.\ Matsumoto. {\em Presentations of subshifts and their topological conjugacy invariants}. Doc.\ Math.\ {\bf 4} (1999), 285-340.

\bibitem{m100} K.\ Matsumoto. {\em $C^*$-algebras associated with presentations of subshifts}. Doc.\ Math.\ {\bf 7} (2002), 1-30.

\bibitem{ma2} K.\ Matsumoto. {\em $C^*$-algebras arising from Dyck systems of topological Markov chains}. preprint, arXiv.math.OA/0607518

\bibitem{mey1}  T.\ Meyerovitch. {\em Tail invariant measures of the Dyck shift}.
preprint, arXiv.math.DS/0406045v2

\bibitem{cbms} I. Raeburn, {\em Graph Algebras}, CBMS Regional Conference Series in Mathematics, vol. 103, Amer. Math. Soc., 2005.


\bibitem{t1} M. Tomforde. {\em A unified approach to Exel-Laca algebras and $C^*$-algebras associated
to graphs}. J. Operator Theory {\bf 50} (2003), 345--368.

\bibitem{t2} M. Tomforde. {\em Simplicity of ultragraph algebras}. Indiana Univ. Math. J. {\bf 52} (2003), 901--926.

\end{thebibliography}
\end{document}